\documentclass{amsproc}

\usepackage{amssymb, amsmath, amscd, latexsym, mathrsfs}
\usepackage{amsfonts}
\usepackage{appendix}
\usepackage[cmtip,all]{xy}

\newtheorem{thm}{Theorem}[subsection]
\newtheorem{lem}[thm]{Lemma}

\newtheorem{prop}[thm]{Proposition}

\theoremstyle{definition}
\newtheorem{defn}[thm]{Definition}
\newtheorem{expl}[thm]{Example}

\theoremstyle{remark}
\newtheorem{rem}[thm]{Remark}
\numberwithin{equation}{subsection}  

\newcommand{\lra}{\longrightarrow}
\newcommand{\co}{\colon\!}

\newcommand{\smin}{\smallsetminus}

\newcommand{\id}{\textup{id}}

\newcommand{\holim}{\textup{holim}}

\newcommand{\hofiber}{\textup{hofiber}}
\newcommand{\mor}{\textup{mor}}

\newcommand{\ob}{\textup{ob}}
\newcommand{\map}{\textup{map}}
\newcommand{\rmap}{\texttt{R}\textup{map}}  
\newcommand{\rsec}{\texttt{R}\Gamma} 

\newcommand{\homeo}{\textup{homeo}}
\newcommand{\aut}{\textup{aut}}

\newcommand{\emb}{\textup{emb}}
\newcommand{\config}{\mathsf{con}}

\newcommand{\fin}{\mathsf{Fin}} 
\newcommand{\finplus}{{\fin_*}}
\newcommand{\holink}{\textup{holink}}
\newcommand{\tw}{\textup{tw}}
\newcommand{\heq}{\textup{he}}

\newcommand{\sC}{\mathscr C}

\newcommand{\sA}{\mathcal A}
\newcommand{\sB}{\mathcal B}
\newcommand{\sD}{\mathcal D}
\newcommand{\sE}{\mathcal E}

\newcommand{\sP}{\mathcal P}
\newcommand{\sQ}{\mathcal Q}

\newcommand{\sU}{\mathcal U}
\newcommand{\op}{\textup{op}}

\newcommand{\haut}{\textup{haut}}
\newcommand{\shift}{E^{\sigma}}
\newcommand{\shiftt}{E^{\tau}}
\newcommand{\col}{\textup{col}}

\newcommand{\NN}{\mathbb N}
\newcommand{\RR}{\mathbb R}

\newcommand{\uli}{\underline}
\newcommand{\Tot}{\textup{Tot}}

\newcommand{\holimsub}[1]{\begin{array}[t]{cc} \textup{holim} \\ [-1mm]
\scriptstyle{#1} \end{array}}

\newcommand{\hocolimsub}[1]{\begin{array}[t]{cc} \textup{hocolim} \\
[-1.7mm] \scriptstyle{#1} \end{array}}

\begin{document}

\title{Configuration categories and homotopy automorphisms}
\author{Michael S. Weiss}%
\address{Math.~Institut, Universit\"at M\"{u}nster, 48149 M\"{u}nster, Einsteinstrasse 62, Germany}%
\email{m.weiss@uni-muenster.de}

\subjclass[2010]{Primary 57R19, 55P65}

\begin{abstract} Let $M$ be a smooth compact manifold with boundary. Under some geometric conditions on $M$, a homotopical
model for the pair $(M,\partial M)$ can be recovered from the configuration category of $M\smin\partial M$.
The grouplike monoid of derived homotopy automorphisms of the configuration category of $M\smin\partial M$
then acts on the homotopical model of $(M,\partial M)$. That action is compatible with a better known homotopical action of the
homeomorphism group of $M\smin\partial M$ on $(M,\partial M)$.
\end{abstract}

\thanks{This project was supported by the Humboldt foundation through a Humboldt professorship, 2012-2017.}
\maketitle


\section{Introduction} \label{sec-intro}
The term \emph{configuration category} of a topological manifold $M$ has a number of interpretations \cite{BoavidaWeissLong}.
In one of them, which is compelling because it makes a direct connection with configuration spaces, it is a category
enriched in topological spaces such that the object space is
\[  \coprod_{k\ge 0} \emb(\uli k\,,M) \]
where $\uli k=\{1,2,\dots,k\}$. A \emph{morphism} from an embedding $f\co \uli k\to M$ to an embedding
$g\co \uli\ell \to M$ consists of a map
\[ v\co \uli k\to \uli\ell~, \]
not necessarily injective, and a (Moore) homotopy $(\gamma_t)_{t\in[0,a]}$ from $f$ to $gv$
which satisfies the stickiness condition: if $\gamma_s(x)=\gamma_s(y)$ for some $s\in [0,a]$ and some $x,y\in \uli k$~, then
$\gamma_t(x)=\gamma_t(y)$ for all $t\in [s,a]$.

For homotopy theoretic purposes it is wise to replace the topological category by its topological nerve, which is a
simplicial space. Therefore $\config(M)$ is strictly speaking a simplicial space, and moreover a \emph{Segal space};
this is one way to say it has the homotopical properties that one expects from the nerve of a well-behaved topological category.
Other ``models'' of the configuration
category of $M$ described in \cite{BoavidaWeissLong} are other simplicial spaces which are degreewise weakly equivalent
to this incarnation. In all these models, $\config(M)$ is a simplicial space \emph{over} the nerve $N\fin$ of $\fin$,
where $\fin$ is the small category whose objects are the finite sets $\uli k$ for $k\ge 0$, and all maps between these
sets as morphisms. As such, $\config(M)$ is a fiberwise complete Segal space over $N\fin$.

A number of people, but especially Bill Dwyer and Ricardo Andrade, have asked whether $\config(M)$ as a homotopical
construct is a good substitute for the topological type of $M$. (In the Dwyer formulation the question probably did
not exactly mention $\config(M)$ but something closely related from the world of operads; and perhaps it was about
$\RR^m$ rather than a general $M$.) More precisely, there is an inclusion map of topological (grouplike) monoids
\begin{equation}  \label{eqn-Dwyer} \homeo(M)  \lra \haut_{N\fin}(\config(M))  \end{equation}
where $\haut_{N\fin}(\config(M))$ denotes the grouplike topological monoid of right derived homotopy automorphisms of
$\config(M)$, as a simplicial space over $N\fin$. One may wonder whether this map is a homotopy equivalence, or a good
approximation in a weaker sense.

For example, \cite[\S10]{BoavidaWeissLong} gives a positive answer in a special case of an analogous question for
manifolds with boundary. Namely, the space of homeomorphisms of a disk $D^m$ relative to the boundary is contractible by
the Alexander trick. The space of homotopy automorphisms of the corresponding configuration category (relative to an
appropriate boundary configuration category and over the nerve of $\finplus$, the appropriate enlargement of $\fin$) is
also contractible.

\medskip
This paper is a continuation of \cite{WeissSen1}. The main point is a translation of some of the results
in \cite{WeissSen1} into a more homotopical language, specifically, the language of Segal spaces. 
After the translation, we have a positive answer to a weak variant of the Dwyer-Andrade question. Here is a
description of that answer (and the weaker question) in a simple case. Let $M$ be a compact smooth manifold with boundary. Let
$M_-= M\smin\partial M$. The boundary $\partial M$ can be recovered from $M_-$ in a homotopical sense as the homotopy link
of the base point in $M/\partial M\cong M_-\cup\infty$. Therefore it is allowed to say that the homeomorphism group
$\homeo(M_-)$ acts on the pair $(M,\partial M)$ by homotopy automorphisms, in the $A_\infty$ sense.

\smallskip
Theorem~\ref{thm-factabs}. \emph{If $M$ is the total space of a smooth disk bundle with fibers of dimension
$\ge 3$ on a smooth closed manifold $L$, then this $A_\infty$ action by homotopy automorphisms
of $\homeo(M_-)$ on the pair $(M,\partial M)$ extends to
an $A_\infty$ action of $\haut_{N\fin}(\config(M_-))$ on $(M,\partial M)$.} 

\smallskip
There is also a more complicated relative version for a smooth compact manifold with boundary and corners.
For that, suppose that $M$ is smooth compact and $\partial M$ is the union of smooth codimension zero
submanifolds $\partial_0M$ and $\partial_1M$ such that $\partial\partial_1M=\partial_0M\cap \partial_1M=
\partial\partial_0M$. The corner set is $\partial_0M\cap \partial_1M$.
We (re)define $M_-$ as $M\smin\partial_1M$ and write $\homeo(M_-;\partial_0M_-)$
for the group of homeomorphisms of $M_-$ which fix the boundary $\partial_0M_-=\partial_0M\smin \partial_1M$
of $M_-$ pointwise.
Let $U$ be an open collar neighborhood
of $\partial_0M$ in $M$ (so that $U\cap\partial_1M$ is an open collar neighborhood of $\partial\partial_1M$ in
$\partial_1M$) and let $U_-:=U\smin\partial_1M$.

\smallskip
Theorem~\ref{thm-factrel}. \emph{If $M$ is a smooth thickening \cite[\S3.2]{WeissSen1} of a neatly embedded compact
smooth submanifold
of $M_-$ of codimension $\ge 3$, then the $A_\infty$ action by homotopy automorphisms of $\homeo(M_-;\partial_0M_-)$ on
the pair $(M,\partial M)$, relative to $\partial_0M$, extends to an $A_\infty$ action of
$\haut_{N\fin}(\config(M_-);\config(U_-))$ on $(M,\partial M)$, still relative to $\partial_0M$.}

\smallskip
Note that theorem~\ref{thm-factabs} is a special case of theorem~\ref{thm-factrel}, the case where
$\partial_0M$ is empty. --- Results of this type are used in \cite{WeissDalian}, so that they have
a \emph{raison d'\^etre} although what they tell us about~(\ref{eqn-Dwyer}) may seem unexciting.
They are likely to be generally useful in manifold calculus applied to spaces of smooth embeddings. They come
with estimates saying that if we are happy to replace $(M,\partial M)$ by its $k$-type (or to kill at least some homotopy
groups above level $k$), then we can restrict attention to configurations in $M$ of cardinality $\le f(k)$, where $f$ is
a fairly uncomplicated function of the variable $k$. 

These estimates become trivial when $\partial M=\emptyset$. In that case it is clear that a homotopy automorphism of
the category of configurations of cardinality $\le 1$ in $M$ determines a homotopy automorphism of $M$. (So
we can take $f(k)=1$ for all $k$, in the notation above.) 

\medskip
Here is a short review of the notation used and the type of results proved in \cite{WeissSen1}, covering the simplest cases.
Suppose first that $L$ is a smooth compact submanifold of a smooth manifold $M$.
Both $M$ and $L$ are without boundary and $L$ is equipped with a Riemannian metric. The elements
of $\sP(L)$ are pairs $(S,\rho)$ where $S$ is a finite subset of $L$ and $\rho\co S\to \RR$ is a function with
positive values. There is a condition: for each $s\in S$, the exponential map is defined and regular on the closed disk of
radius $\rho(s)$ in $T_sL$, and the images of these disks under the exponential maps (for each $s\in S\subset L$) are
pairwise disjoint in $L$. Let $V_L(S,\rho)\subset L\subset M$ be the union of the images in $L$ of the corresponding
\emph{open} balls of radius $\rho(s)$ in $T_sL$ under the exponential maps. There are results of the following form: the map
\[  M\smin L \lra \holimsub{(S,\rho)\in \sP(L)} M\smin V_L(S,\rho) \]
induced by the inclusions $M\smin L\to M\smin V_L(S,\rho)$ is a weak equivalence, under some conditions.
That map can also be written in the form
\[ M\smin L \lra \holim~\Phi \]
where $\Phi$ is the contravariant functor $(S,\rho)\mapsto M\smin V_L(S,\rho)$ from $\sP(L)$ to spaces.
The homotopy limit is an enriched variant. \newline Suppose next that $M$ is smooth, compact, \emph{with} boundary and
equipped with a Riemannian metric; no submanifold $L$ is specified. Define $\sP(M\smin \partial M)$ roughly as above in the
case of $L$, so that elements of $\sP(M\smin\partial M)$ are pairs $(S,\rho)$ where $S$ is a finite subset of $M\smin\partial M$
and $\rho$ is a function with positive values on $S$. Again there is a condition: for each $s\in S$, the exponential map
is defined and regular on the closed disk of radius $\rho(s)$ in $T_s(M\smin\partial M)$, and the images of these disks
under the exponential maps (for each $s\in S$) are pairwise disjoint in $M\smin\partial M$. For
$(S,\rho)\in \sP(M\smin\partial M)$ let
$V(S,\rho)\subset M\smin\partial M$ be the union of the images of the corresponding \emph{open} balls of radius $\rho(s)$ in
$T_x(M\smin\partial M)$ under the exponential maps. There are results of the following form: the map
\[  \partial M \lra \holimsub{(S,\rho)\in \sP(M)} M\smin V(S,\rho) \]
induced by the inclusions $\partial M\to M\smin V(S,\rho)$
is a weak equivalence, under some conditions. That map can also be written in the form
\[ \partial M \lra \holim~\Psi \]
where $\Psi(S,\rho)= M\smin V(S,\rho)$ for $(S,\rho)\in \sP(M\smin\partial M)$.
The homotopy limit is an enriched variant.

\medskip
From the point of view that we take here, there are several aspects to these results 
in \cite{WeissSen1} which call for a translation. Clearly we want to have a more flexible
replacement for $\sP(M\smin\partial M)$, one which does not require a smooth structure on $M$.
This replacement is the configuration category of $M\smin\partial M$ in various guises.
We view it as a Segal space; more about all that in section~\ref{sec-langcon}. But that is not the
demanding part of the translation. A more important aspect is that \cite{WeissSen1} gives a homotopical description of $\partial M$
(under conditions) in terms of $\sP(M\smin\partial M)$, or the configuration category of $M$, \emph{and} a
contravariant functor $\Psi$ on $\sP(M\smin\partial M)$; but here we need a homotopical description of $\partial M$
in terms of the configuration category of $M$ \emph{only}. Therefore it is on us to explain
why the functor $\Psi$ is intrinsic to the configuration category. The explanation is not difficult,
but it can be given in very general terms and therefore it will be given in very general terms.
This is going to happen mainly in section~\ref{sec-shifty}.

\section{Segal spaces and configuration categories} \label{sec-langcon}  
The translation promised in the introduction uses Rezk's concept of
Segal space and the associated framework \cite{Rezk} in which topological categories
can be viewed as objects of a model category. It is not a great challenge to recast the topological posets
$\sP(L)$, $\sP(M\smin\partial M)$ etc.~as
complete Segal spaces. Indeed their topological nerves are already complete Segal spaces; but we are going to tinker with
the definitions in order to make a better connection with \cite{BoavidaWeissLong}. 

\subsection{Simplicial spaces, Segal spaces and complete Segal spaces}
In the first few sections of \cite{BoavidaWeissLong} there are definitions of \emph{Segal space}
and \emph{complete Segal space} with some examples and failing candidates. The central example
there is the configuration category of a manifold $M$. This comes in many guises, well-known and less well-known.
There is a variant where $M$ has empty boundary, and a more complicated variant in the case where $M$ has
nonempty boundary. Here we only give a brief review of the definitions and examples.

\bigskip
The nerve construction turns small categories into simplicial sets and small topological categories
into simplicial spaces. It was Graeme Segal \cite{Segal}, \cite{Adams} who promoted the idea that the nerves and their
homotopical properties are more important than the categories themselves. In that spirit Rezk \cite{Rezk}
gave the following definition. For $r\ge 0$ and $i=1,2,\dots,r$ let $u_i\co[1]\to [r]$ be the
monotone injection defined by $u_i(0)=i-1$ and $u_i(1)=i$.
A \emph{Segal space} is a simplicial space $X$ such that for each $r\ge 2$ the map
$(u_1^*,u_2^*,\dots,u_r^*)$
from $X_r$ to the homotopy inverse limit of the diagram
\[
\xymatrix{
X_1\ar[r]^-{d_0} & X_0 & \ar[l]_{d_1} X_1 \ar[r]^-{d_0} & \cdots & \cdots\ar[r]^-{d_0}& X_0 & \ar[l]_-{d_1}  X_1\,,
}
\]
with $r$ copies of $X_1$~, is a weak homotopy equivalence. (In the case where $X_0$ is weakly contractible
this simplifies to the condition
that $(u_1^*,u_2^*,\dots,u_r^*)$ as a map from $X_r$ to $(X_1)^r$ be a weak homotopy equivalence for all $r\ge 2$. This
constitutes Segal's definition or interpretation of what it means for the space $X_1$ to have the structure of
an $A_\infty$ topological monoid with unit.) In particular the nerve of
any small category is a Segal space $X$ which also happens to be a simplicial \emph{set}. Another
important type of example: if $\sC$ is a topological category
(category object in the category of topological spaces), and if one of the maps \emph{source}, \emph{target}
from the space $\mor(\sC)$ to the space $\ob(\sC)$ is a Serre fibration, then the nerve $N\sC$ is also a Segal space.

\medskip
\emph{Terminology and notation.} Let $X$ be a Segal space. We view $X_0$ as the space of objects, $X_1$ is the space of morphisms
and $d_0,d_1\co X_1\to X_0$ as the operators \emph{source} and \emph{target}. For elements $a,b\in X_0$
the space $\mor_X^h(a,b)$ is the homotopy fiber of $(d_0,d_1)\co X_1\to X_0\times X_0$
over $(a,b)$.  

\medskip
Equivalences between small categories are not always reflected in degreewise
weak equivalences between their nerves. Indeed if $\sC$ is equivalent to $\mathscr D$, then there
is no strong reason why
$N_0\sC=\ob(\sC)$ should be weakly equivalent to $N_0\mathscr D=\ob(\mathscr D)$. To deal with this,
Rezk introduced the concept of
\emph{Dwyer-Kan} equivalence between Segal spaces as an analogue of the classical concept of equivalence
between categories (see definition~\ref{defn-DwyerKan} just below) and a related concept of completeness. A Segal space $Y$ is \emph{complete} if
the map $d_0\co Y_1\to Y_0$ alias \emph{source} restricts to a weak equivalence from $Y_1^\heq$ to $Y_0$~, where $Y_1^\heq$
is the union of the homotopy invertible path components of $Y_1$. He showed that for every Segal space $X$ there exists
a Segal space $Y$ and a simplicial map $X\to Y$ which is a Dwyer-Kan equivalence
and where $Y$ is complete.
Moreover a Dwyer-Kan equivalence between complete Segal spaces is automatically a degreewise weak equivalence (between simplicial
spaces). \emph{Example}: a discrete group $G$ determines a category with one object whose endomorphism monoid
is the group $G$. The nerve of that category is a Segal space $X$, but it is not complete
unless $G$ is the trivial group. The Rezk completion $Y$ of $X$ has the form of a constant simplicial space,
$Y_r=BG$ for all $r\ge 0$.

\begin{defn} \label{defn-DwyerKan} A map $f\co X\to Y$ between Segal spaces is a \emph{Dwyer-Kan equivalence} if
\begin{itemize}
\item[-] for all $a,b\in X_0$ the induced map $\mor^h_X(a,b)\to \mor^h_Y(f(a),f(b))$ is a weak equivalence;
\item[-] for every $b\in Y_0$ there exist $a\in X_0$ and $c\in Y_1^\heq$ such that
$d_1c$ is in the path component of $b$ and $d_0c$ is in the path component of $f(a)$.
\end{itemize}
\end{defn}

\medskip
The nerves of the topological posets $\sP(L)$ and $\sP(M\smin\partial M)$ defined in section~\ref{sec-intro}
are examples of complete Segal spaces. In \cite{BoavidaWeissLong} we use slightly
different editions denoted $\config(L)$, $\config(M)$, $\config(M\smin\partial M)$ etc., for mostly bureaucratic reasons. One
definition of the simplicial space $X=\config(L)$ for a smooth Riemannian manifold $L$ with \emph{empty} boundary is as follows.
An object, alias element of $X_0$~, is an element $(S,\rho)$ of $\sP(L)$ together with a total ordering of $S$.
There is at most one morphism between any two objects, and this happens if and only
if $(S,\rho)\le (T\,\sigma)$ holds for the underlying elements $(S,\rho)$ and $(T,\sigma)$ of $\sP(L)$.
Thus $X_0$ is a covering space of the space of objects of
$\sP(L)$, so that the fiber over an element $(S,\rho)$ in $\sP(L)$ is the set of total orderings of $S$;
and $X_1$ is a covering space of the space of morphisms in $\sP(L)$, so that the fiber over $(S,\rho)\le (T,\sigma)$
is the product of the set of total orderings of $S$ with the set of total orderings of $T$. In this way, $X_0$ and $X_1$
form object space and morphism space of a topological category (category object in the category of topological spaces).
The nerve of that is $\config(L)$, a simplicial space; more precisely it is called the \emph{Riemannian model} of
$\config(L)$ in \cite{BoavidaWeissLong}. It is a Segal space
(e.g.~because $d_1=$\emph{target} from $X_1$ to $X_0$ is a Serre fibration),  but not a complete Segal space
except in a few cases of little interest. We recover
the loss by making two observations.
\begin{itemize}
\item The forgetful functor $\config(L)\to \sP(L)$ is a Dwyer-Kan equivalence.
\item $\config(L)$ is a \emph{fiberwise complete Segal space} over the nerve of $\fin$
(explanation follows).
\end{itemize}
Here $\fin$ is the small category whose objects are the finite sets
$\uli n=\{1,2,\dots,n\}$ for $n\ge 0$ with all maps between these sets as morphisms. There is an obvious
forgetful functor from $\config(M)$ in the above definition to $N\fin$, the nerve of $\fin$.
By saying that $\config(L)\to N\fin$ is a fiberwise complete Segal space we mean that the resulting
commutative square
\[
\xymatrix@R=12pt{ X_1^\heq \ar[r]^-{d_0} \ar[d] & \ar[d] X_0 \\
Y_1^\heq  \ar[r]^-{d_0} & Y_0
}
\]
(where $X=\config(L)$ and $Y=N\fin$) is a weak homotopy pullback square. \newline
In the case where $L$ has nonempty boundary, there is a more complicated definition of $\sP(L)$
and a related definition of $\config(L)$. The elements
of $\sP(L)$ are pairs $(S,\rho)$ where $S$ is a finite subset of $L\smin\partial L$ and $\rho$ is a
function from $S\sqcup \partial L$ to the positive reals, locally constant
on $\partial L$ and subject to a few more conditions.
\begin{itemize}
\item[-] For each $s\in S$, the exponential map $\exp_s$ at $s$ is defined and regular on the disk
of radius $\rho(s)$ about the origin in $T_sL$\,.
\item[-] The (boundary-normal) exponential map is defined and regular on the set of all tangent vectors $v\in T_zL$ where $z\in \partial L$,
where the vector $v$ is inward perpendicular to $T_z\partial L$ and $|v|\le \rho(z)$.
\item[-] The images in $L$ of these disks and the image
of this band under the exponential map(s) are pairwise disjoint.
\end{itemize}
For a pair $(S,\rho)$ satisfying these conditions, let $V(S,\rho)\subset L$ be the union of the
open balls of radius $\rho(s)$ about elements $s\in S$ and the open collar on $\partial L$ determined the normal distance
function $\rho|_{\partial L}$. (Sometimes we write $V_L(S,\rho)$ instead of $V(S,\rho)$, for example if $L$ comes as a
smooth submanifold of another smooth manifold $M$.) In this case
there is a reference map $\config(L)\to N\finplus$, where $\finplus$ is the following category.
Objects are the finite sets $[n]=\{0,1,\dots,n\}$ for $n\ge 0$ which we view as sets with a base point $0$,
and the morphisms are all base-point preserving maps between these. Both $\config(L)$ and $N\finplus$
are Segal spaces and the reference map $\config(L)\to N\finplus$ makes $\config(L)$ into a
fiberwise complete Segal space over $N\finplus$.

\bigskip
There are many alternative descriptions of $\config(L)$ in \cite{BoavidaWeissLong}; each of these can be related
to the above definition by a chain of degreewise equivalences over $N\fin$ or over $N\finplus$~, as appropriate.
(Actually $L$ is typically called $M$ in \cite{BoavidaWeissLong} and now we shall adopt that habit.)
One of them, the \emph{particle model}, deserves to be mentioned here because it does not require a Riemannian
metric, or even a smooth structure, and has very good naturality properties. It is probably due to \cite{Andrade}.
Suppose to begin with that $M$ is a topological manifold with empty boundary.
Let $k\in\NN$. The space of maps from $\uli k$ to $M$
comes with an obvious stratification. There is one stratum for each equivalence relation $\eta$ on $\uli k$~. The
points of that stratum are precisely the maps $\uli k\to M$ which can be factorized as projection from $\uli k$ to
$\uli k\,/\eta$ followed by an injection of $\uli k\,/\eta$ into $M$. \newline
Now we construct a topological category whose object space is
\begin{equation} \label{eqn-appfact1}  \coprod_{k\ge 0} \emb(\uli k\,,M). \end{equation}
By a morphism from $f\in \emb(\uli k\,,M)$ to $g\in\emb(\uli\ell\,,M)$ we mean a pair consisting of a map $v\co \uli k\to \uli\ell$
and a \emph{reverse exit path} $\gamma$ from $f$ to $gv$ in the stratified space of \emph{all} maps from $\uli k$ to $M$.
(In more detail: $\gamma$ is a path $[0,a]\to \map(\uli k,M)$ for some $a\ge 0$, and the reverse exit path property means
that if $\gamma_t(x)=\gamma_t(y)$ for some $t\in [0,a]$ and $x,y\in S$, then $\gamma_s(x)=\gamma_s(y)$
for all $s\in[t,a]$. Note that $f$ is injective but $gv$ need not be injective since $v$ is not required
to be injective.) The space
of all morphisms is therefore a coproduct with one summand for each morphism $v\co \uli k\to \uli \ell$ in $\fin$,
where that summand consists of triples $(f,g,\gamma)$ as above: $f\in\emb(\uli k\,,M)$,
$g\in\emb(\uli\ell\,,M)$ and $\gamma$ a reverse exit path from $gv$ to $f$. Composition of morphisms
is obvious. The nerve of this category is a fiberwise complete Segal space
over $N\fin$.

In the case of a topological manifold $M$ with nonempty boundary, the definition of
$\config(M)$ along similar lines is slightly more complicated, but we need it.
The space of maps from $\uli k$ to $M$ comes with a stratification.
There is one stratum for each pair $(S,\eta)$ where $S\subset\uli k$ and $\eta$ is an equivalence relation
on $\uli k$ such that $S$ is a union of equivalence classes. The
points of that stratum are the maps $\uli k\to M$ taking $S$ to $\partial M$ and the complement
of $S$ to $M\smin\partial M$, and which can be factored as projection from $\uli k$ to $\uli k/\eta$
followed by an injection of $\uli k/\eta$ into $M$.
Now we construct a topological category whose object space is
\begin{equation}  \coprod_{k\ge 0} \emb(\uli k\,,M\smin\partial M). \end{equation}
A morphism from $f\in \emb(\uli k\,,M\smin\partial M)$ to $g\in\emb(\uli\ell\,,M\smin\partial M)$ is a
pair consisting of a morphism $v\co[k]\to[\ell]$ in $\finplus$
and a Moore path $\gamma=(\gamma_t)_{t\in[0,a]}$ in $\map(\uli k\,,M)$ which is a reverse exit path. It is required to
satisfy $\gamma_0=f$ and $\gamma_a(x)=g(v(x))$ if $v(x)\in\uli\ell$\,, but $\gamma_a(x)\in \partial M$ if $v(x)=0$.
Composition of morphisms is almost obvious. The nerve of this topological category is a fiberwise complete Segal space
over $N\finplus$ which we can regard as an alternative definition or description of $\config(M)$.

\subsection{Functors as maps between simplicial spaces}
A \emph{map} between complete Segal spaces $X$ and $Y$ is a simplicial map $f\co X\to Y$. Such an $f$ can also be
regarded as a functor from $X$ to $Y$. The map $f$ is considered to be a \emph{weak equivalence} if each $f_r\co X_r\to Y_r$
is a weak equivalence of spaces. For many purposes it is useful to have a notion of \emph{space} of maps from $X$
to $Y$ which is functorial in the two variables $X$ and $Y$ and takes weak equivalences to weak equivalences.
Such a concept exists and is called the \emph{derived} space of maps from $X$ to $Y$, and denoted
\[  \rmap(X,Y)~. \]
To define this we do not need to know or assume that $X$ and $Y$ are complete Segal spaces. It suffices to know
that they are simplicial spaces. It suffices to have a decision as to which simplicial maps between simplicial
spaces are to be called \emph{weak equivalences} (namely, those which are degreewise weak equivalences of spaces). \newline
Therefore we switch briefly to the general setting where $X$ and $Y$ are contravariant functors from a
(small, discrete) category $\sC$ to the category of spaces. (The example to keep in mind is $\sC=\Delta$.)
Let $\mathscr D$ be the category of such functors from $\sC$ to spaces, where morphisms alias \emph{maps}
are natural transformation. A map $f\co X\to Y$ in $\mathscr D$ is a weak equivalence if $f_c\co X(c)\to Y(c)$
is a weak equivalence of spaces for each object $c$ in $\sC$. It is straightforward to define a space
$\map(X,Y)$, for example as the geometric realization of the simplicial set where a $k$-simplex is a map from
$X\times\Delta^k$ to $Y$, where $\Delta^k$ is the geometric $k$-simplex (a space).
We look for a definition of $\rmap(X,Y)$, the derived mapping space. There are two
well-known options.
\begin{itemize}
\item
Dwyer and Kan \cite{DwyKa2} have a definition of $\rmap(X,Y)$ in an extremely general setting where
$X$ and $Y$ are objects in a category $\mathscr D$ with a subcategory of so-called weak equivalences,
subject to some mild conditions. Their definition of $\rmap(X,Y)$ is big in the sense that it can be a simplicial class
rather than a simplicial set if $\mathscr D$ is not small.
\item For the category $\mathscr D$ (as defined above, category of contravariant functors
from $\sC$ to spaces) we have a subcategory of \emph{weak equivalences} and a preferred action of the category
of simplicial sets on $\mathscr D$, given by the ordinary degreewise product of simplicial sets
with simplicial spaces. There are a few standard ways to enhance these data to the structure
of a Quillen simplicial model category \cite{Hovey}, \cite{Hirschhorn}. (For us the preferred choice is the one where a map
between simplicial spaces is considered to be a fibration if it is a degreewise Serre fibration.)
Then  $\rmap(X,Y)$ can be defined as $\map(X^\flat,Y^\sharp)$ where $X^\flat$
is a cofibrant replacement of $X$ and $Y^\sharp$ is a fibrant replacement of $Y$.
To achieve strict functoriality one should use functorial replacements, so that $X\mapsto X^\flat$
is a functor with a natural transformation to the identity by weak equivalences, and $Y\mapsto Y^\sharp$ is a functor
with a natural transformation from the identity by weak equivalences.
\end{itemize}
It is a special case of a result in \cite{DwyKa2} that these two definitions of $\rmap(X,Y)$ agree up to a chain of weak
equivalences. (This has the consequence that $\rmap(X,Y)$ according to the second definition is largely
independent of the choices required there.)

\bigskip
Returning to simplicial spaces $X$ and $Y$, we conclude that we have a few good definitions
of a derived mapping space $\rmap(X,Y)$, since a simplicial space is a contravariant functor from $\Delta$ to spaces.
More generally, suppose that $Z$ is a fixed simplicial set, and let $X,Y$ be simplicial spaces over $Z$,
that is to say, simplicial spaces equipped with reference maps $p_X$ and $p_Y$ to $Z$, respectively.
By $\rmap_Z(X,Y)$ we mean the fiber of the map $\rmap(X,Y)\to \map(X,Z)$ given by composition with $p_Y$
over the point determined by $p_X$. Perhaps it is worth pointing out that $\map(X,Z)$ is a set, alias discrete
space. We are also using the fact that $p_Y\co Y\to Z$ extends uniquely to a map $Y^\sharp\to Z$.

\subsection{Grothendieck construction and homotopy inverse limits}  
Let $\sC$ be a category, small and discrete for simplicity, and let $F$ be a functor from $\sC^\op$ to
spaces. The \emph{Grothendieck construction} of $F$ is a new category $\smallint F$ internal to the category
of spaces. It comes with a forgetful functor $\smallint F \lra \sC$. The space of objects of $\smallint F$
is
\[  \coprod_{c\textup{ in }\sC} F(c) \]
so that we may describe objects as pairs $(c,x)$ where $c$ is an object of $\sC$ and $x\in F(c)$.
The space of morphisms in $\smallint F$ is
\[  \coprod_{g\co c\to d} F(d) \]
where $g$ runs through all the morphisms in $\sC$. We may describe morphisms in $\smallint F$ as pairs
$(f,x)$ where $f\co c\to d$ is in $\sC$ and $x\in F(d)$. The source of a morphism $(f\co c\to d,x)$ in $\smallint F$
is $(c,g^*x)$; the target is $(d,x)$. The forgetful functor to $\sC$ is given on objects by $(c,x)\mapsto c$ and on
morphisms by $(g,x)\mapsto g$.

The space of sections of the forgetful functor $\smallint F\to \sC$ is identifed with the space $\lim~F$
(inverse limit of $F$). This can also be formulated with nerves. The forgetful functor induces a map of nerves
\[  N\smallint F\lra N\sC \]
and the space of sections of that (in the category of simplicial spaces) is identified with the space $\lim~F$.

Furthermore, a map $p\co X\to N\smallint F$ of simplicial spaces is, up to isomorphism over $N\smallint F$,
the projection associated with a Grothendieck construction $\smallint F$ (for $F$ from $\sC^\op$ to spaces)
if and only if
\begin{itemize}
\item[(i)] $X$ is a nerve, that is to say, for each $r\ge 2$ the map
$(u_1^*,u_2^*,\dots,u_r^*)$ from $X_r$ to the inverse limit of the diagram
\[
\xymatrix{
X_1\ar[r]^-{d_0} & X_0 & \ar[l]_{d_1} X_1 \ar[r]^-{d_0} & \cdots & \cdots\ar[r]^-{d_0}& X_0 & \ar[l]_-{d_1}  X_1\,,
}
\]
with $r$ copies of $X_1$~, is a homeomorphism;
\item[(ii)] the commutative diagram
\[
\xymatrix{
X_0 \ar[d]^-p & \ar[l]_{d_1} X_1 \ar[d]^-p \\
(N\sC)_0 & \ar[l]_{d_1} (N\sC)_1
}
\]
is a pullback square of topological spaces.
\end{itemize}

\medskip
We take this characterization of Grothendieck constructions as a model for a homotopical (derived) variant
of the Grothendieck construction.

\begin{defn} \label{defn-rightfib} Let $X$ and $Y$ be Segal spaces. A map $p\co X\to Y$ (in the category of simplicial spaces)
is a \emph{right fibration} if
the commutative diagram
\[
\xymatrix{
X_0 \ar[d]^-p & \ar[l]_{d_1} X_1 \ar[d]^-p \\
Y_0 & \ar[l]_{d_1} Y_1
}
\]
is (weakly) homotopy cartesian.
\end{defn}
This definition is taken from \cite{BoavidaAlpine}. It has a better known counterpart in the setting
of quasi-categories alias $\infty$-categories, as explained in \cite{BoavidaAlpine}.

\smallskip
In the situation of~\ref{defn-rightfib} we may loosely think of $X$ as the Grothendieck construction of a contravariant functor from $Y$
to spaces whose value at an object $b\in Y_0$ is the homotopy fiber of $p\co X_0\to Y_0$ over $b$.
A morphism $g\in Y_1$ with source $a=d_0g$ and target $b=d_1g$ determines a diagram
\[
\xymatrix{ \hofiber_b[X_0\to Y_0] & \ar[l]_{d_1}^\simeq \hofiber_g[X_1\to Y_1] \ar[r]^-{d_0} & \hofiber_a[X_0\to Y_0]
}
\]
which can be viewed as a map in the derived sense
from $\hofiber_b[X_0\to Y_0]$ to $\hofiber_a[X_0\to Y_0]$.

\begin{rem} In the notation of~\ref{defn-rightfib}, if we think of $p\co X\to Y$ as the projection associated with
a Grothendieck construction corresponding to an elusive functor (from an elusive category with nerve $Y$ to spaces),
then we must think of the derived section space of $p$, denoted $\rsec(p)$, as the homotopy limit of the elusive functor.
\end{rem}

To spell out what $\rsec(p)$ is, we select a simplicial model category structure on the category of simplical spaces where the weak equivalences are the levelwise weak homotopy equivalences. Choose a factorization
\[
\xymatrix@C=16pt{
X \ar[r] & X^\sharp \ar[r]^{p^\sharp} & Y
}
\]
of $p$ where the first arrow is a weak equivalence and $p^\sharp$ is a fibration (in the
selected model category structure).
Choose a weak equivalence $Y^\flat \to Y$ where $Y^\flat$ is cofibrant.
The derived section space of $p$ is the space (simplicial set) of
lifts as in the diagram
\[
\xymatrix{
 & X^\sharp \ar[d]^-{p^\sharp} & \ar[l]_\simeq X \ar[dl]^-p \\
Y^\flat \ar@{..>} [ur] \ar[r]^-\simeq & Y
}
\]

\section{The shift construction} \label{sec-shifty}
What we are after in this section is a
description of functors such as the functor $\Psi$ in section~\ref{sec-intro}, and their homotopy inverse limits,
in terms of not much more than the source category. In the setting of section~\ref{sec-intro} the source
category would be $\sP(M\smin\partial M)$, but it is better for us to use the variant $\config(M\smin\partial M)$ with
the reference map from there to the nerve of $\fin$.

\subsection{The shift construction as an example of a right fibration} \label{subsec-shifty}
Let $X$ be a Segal space and let $A$ be any simplicial space. Rezk has a definition of an
\emph{internal} mapping object $X^A$ which is as follows (in a possibly simplified form which I hope is
good enough here). Put
\[ (X^A)_n:= \rmap(\Delta[n]\times A,X) \]
where $\Delta[n]$ is the simplicial set freely generated by one element in  degree $n$, so that
the geometric realization $|\Delta[n]\,|$ is the standard geometric $n$-simplex $\Delta^n$.

\begin{prop} \emph{(Rezk.)} If $X$ is a Segal space, then $X^A$ is a Segal space; if $X$ is a complete
Segal space, then $X^A$ is a complete Segal space. \qed
\end{prop}

\begin{expl} If $X$ and $A$ are both complete Segal spaces, then $X^A$ should be viewed as the \emph{category}
of functors from $A$ to $X$. In particular $(X^A)_0\cong \rmap(A,X)$ should be viewed
as the (derived) space of functors from $A$ to $X$ and $(X^A)_n=\rmap(\Delta[n]\times A,X)$
should be viewed as the (derived) space of constellations
\[ G_0 \leftarrow G_1 \leftarrow \cdots \leftarrow G_n \]
where $G_0,G_1,\dots,G_n$ are functors from $A$ to $X$
and the arrows connecting them are natural transformations. (In particular, if $A$ is the simplicial
space which has a single point in each degree, then we recover the idea that $(X^A)_n\simeq X_n$
is the derived space of functors from $[n]^\op$ to $X$.)
\end{expl}

\begin{defn} \label{defn-shiftfunctor} There is a functor $\sigma\co\fin\to \fin$ given by disjoint sum with a
singleton. In more
detail, $\sigma$ is given by $\{1,2,\dots,k\}\mapsto \{1,2,\dots,k,k+1\}$ on objects, and for a morphism
$f\co \uli k\to \uli\ell$ the morphism $\sigma(f)$ is given by $\sigma(f)(x)=f(x)$ for $x\le k$
and $\sigma(f)(k+1)=\ell+1$. The standard inclusions of $\{1,2,\dots,k\}$ in $\{1,2,\dots,k,k+1\}$ define
a natural transformation $u\co \id\to \sigma$ between endofunctors of $\fin$\,. Together, $\sigma$ and $u$
determine a map from $\Delta[1]\times N\fin$ to $N\fin$, or in adjoint form, a
map of Segal spaces
\[ N\fin \to N\fin^{\Delta[1]}. \]
Since $N\fin^{\Delta[1]}$
is a simplicial set and at the same time a Segal space, it is (isomorphic) to the nerve of a small category,
and this is also easy to see directly.
\end{defn}

\begin{defn} (\emph{Shift construction.})
Let $X$ be a fiberwise complete Segal space over $N\fin$. Let $\shift(X)$ be the Segal
space defined by the pullback square of simplicial spaces and simplicial sets
\[
\xymatrix{  \shift(X) \ar[d] \ar[r] &  X^{\Delta[1]} \ar[d] \\
N\fin \ar[r]^-{(\sigma,u)} &  N\fin^{\Delta[1]}_{\rule{0mm}{1mm}}
}
\]
There is a map $\psi_X\co \shift(X)\to X$ over $N\fin$ given by composing $\shift(X)\to X^{\Delta[1]}$ from the defining
pullback square with the map
\[ X^{\Delta[1]}\to X^{\Delta[0]}\cong X \]
determined by the
map $\Delta[0]\to \Delta[1]$ which takes the preferred generator in degree $0$ to $d_0$
of the preferred generator in degree $1$.
\end{defn}

Now we look for additional conditions on $X$ which ensure that $\psi_X$ is a right fibration (definition~\ref{defn-rightfib}).
\emph{Notation and vocabulary}: For an $r$-simplex
\[  C = (\uli k_0 \xleftarrow{f_1} \uli k_1 \xleftarrow{f_2} \uli k_2
\xleftarrow{f_3} \cdots \xleftarrow{f_r} \uli k_r) \]
in $N\fin$, denote by $X(C)$ the part of $X_r$ projecting to $C$
under the reference map $X\to N\fin$. An element of $\uli k_j$~, where $j\in \{1,\dots,r-1\}$,
is \emph{heavy} for the diagram $C$ if its preimage under $f_{j+1}$ has more than one element.

\begin{defn} \label{defn-elementary}
A Segal space $X$ over $N\fin$ is \emph{elementary} if it is fiberwise complete and, for
every $r$-simplex
\[   C = (\uli k_0 \xleftarrow{f_1} \uli k_1 \xleftarrow{f_2} \uli k_2
\xleftarrow{f_3} \cdots \xleftarrow{f_r} \uli k_r) \]
in $N\fin$ and $j\in \{1,\dots,r-1\}$ such that $\uli k_j$ has no heavy elements, the
face operator $d_j\co X_r\to X_{r-1}$ restricts to a weak equivalence $X(C) \to X(d_jC)$.
(Because of the Segal property it is enough to require this for $r=2$, in which case $j=1$ is compulsory.)
\end{defn}

\begin{expl} For a manifold $M$, the Segal space $\config(M)$ over $N\fin$ is elementary.
\end{expl}

\begin{prop} \label{prop-elementaryright} Let $X$ be a Segal space over $N\fin$.
If $X$ is elementary, then $\psi_X\co \shift(X)\to X$ is a right fibration. \end{prop}

\proof Let $f\co \uli k\to \uli\ell$ be a morphism in $\fin$. We need to show that
\begin{equation} \label{eqn-elementaryright1}
\begin{aligned}
\xymatrix@R=18pt{
(E^\sigma X)(\uli\ell)  \ar[d]^-\psi & \ar[l]_-{d_1}  (E^\sigma X)(f) \ar[d]^-\psi \\
X(\uli\ell) & \ar[l]_{d_1}  X(f)
}
\end{aligned}
\end{equation}
is weakly cartesian. ---
Let $A$ be the commutative diagram
\[
\xymatrix{
\uli k \ar@{..>}[dr] \ar[d]_-{u_k} \ar[r]^-f & \uli\ell \ar[d]^-{u_\ell} \\
\sigma(\uli k) \ar[r]_-{\sigma(f)}  & \sigma(\uli\ell) \\
}
\]
in $\fin$. We view $A$ as a map from $\Delta[1]^2$ to $N\fin$ and write $X(A)$ for the part
of $\rmap(\Delta[1]^2,X)$ projecting to $A\co \Delta[1]^2\to \fin$. In the category of simplicial spaces, with levelwise
weak equivalences, $\Delta[1]^2$ is the pushout (and also the homotopy pushout) of two copies of $\Delta[2]$ along a common
edge $\Delta[1]$; this corresponds to the decomposition of $A$ into upper triangle $A'$ , lower
triangle $A''$ and diagonal arrow. Now~(\ref{eqn-elementaryright1}) simplifies to
\begin{equation} \label{eqn-elementaryright2}
\begin{aligned}
\xymatrix{
X(u_\ell)  \ar[d]^-{d_0} & \ar[l]_-{\textup{res.}} X(A) \ar[d]^-{\textup{res.}} \\
X(\uli\ell) & \ar[l]_{d_1}  X(f)
}
\end{aligned}
\end{equation}
Using the elementary property of $X$ and the homotopy pushout decomposition of $\Delta[1]^2$ we can also
replace $X(A)$ forgetfully by the weakly equivalent $X(A')$.
The Segal property of $X$ now implies that~(\ref{eqn-elementaryright2}) is weakly cartesian. \qed

\begin{rem} Let $X$ be an elementary Segal space over $N\fin$. Let
\[ C= (\uli k_0 \leftarrow \uli k_1 \leftarrow \cdots \leftarrow \uli k_r) \]
be an $r$-simplex in $N\fin$. Let $D$ be the $(r+1)$-simplex in $\fin$ obtained from $C$
by adding one arrow, the inclusion of $\uli k_0=\{1,2,\dots,k_0\}$ in $\{1,2,\dots,k_0,k_0+1\}$.
Reasoning as in the proof above leads to the following: there is a forgetful weak equivalence from
$E^\sigma(X)(C)$ to $X(D)$.
\end{rem}

\bigskip
\begin{rem} \label{rem-hautsec} Suppose that $X$ is a fiberwise complete Segal space over $N\fin$.
If $X$ is elementary, then $\rsec(\psi_X)$ is a model for the homotopy inverse limit of the contravariant functor on
$X$ encoded by $\psi_X$ according to proposition~\ref{prop-elementaryright}. Whether $X$ is elementary or not, we
would like to say that homotopy automorphisms of $X$ as a simplicial space over $N\fin$
induce homotopy automorphisms of $\rsec(\psi_X)$. The following is a rather pedestrian justification.
[~The idea is that the classifying space of the grouplike
monoid of homotopy automorphisms of $X$ over $N\fin$ carries a universal simplicial fibration with fibers
weakly equivalent to $X$, and we may apply $\rsec(\psi_{(-)})$ fiberwise to that.~]
Abbreviate $W(X):=\rsec(\psi_X)$. Let $v\co X\to Y$ be a weak equivalence between simplicial spaces
over $N\fin$, both fiberwise complete Segal. This determines a commutative square
of simplicial spaces
\[
\xymatrix{
\shift(X) \ar[d]^-{\psi_X} \ar[r]^-{\shift(v)} & \shift(Y) \ar[d]^-{\psi_Y} \\
X \ar[r]^-v & Y
}
\]
Define $W(v)$ as the simplicial set of triples $(s,t,h)$ where $s$ and $t$ are derived sections of $\psi_X$
and $\psi_Y$ respectively, and $h$ is a homotopy connecting $\shift(v)^\sharp\circ s$ to $t\circ v^\flat$. (More precisely
a $k$-simplex of $W(v)$ is a family of such triples $(s,t,h)$ parameterized by the geometric $k$-simplex $\Delta^k$.)
There are forgetful weak equivalences
\[
\xymatrix{ W(X) & \ar[l] W(v) \ar[r] & W(Y)
}
\]
which are also Kan fibrations. This is already enough
to establish a great deal of naturality for the construction $X\mapsto W(X)$.
Namely, for a fixed $X$, simplicial
space over $N\fin$, let $\sC_X$ be a small subcategory of the category of simplicial spaces
over $N\fin$ (no enrichment assumed here) with the following properties.
\begin{itemize}
\item[-] $X$ is an object of $\sC_X$.
\item[-] Every object of $\sC_X$ is weakly equivalent to $X$ (in the category of simplicial spaces over $N\fin$).
\item[-] If $Y$ belongs to $\sC_X$\,, then
$Y\times\Delta^k$ also belongs to $\sC_X$ for every $k\ge 0$.
\item[-] The morphisms in $\sC_X$ between any two objects of $\sC_X$ are precisely the weak equivalences between these two in the
category of simplicial spaces over $N\fin$.
\item[-] $\sC_X$ is closed under (some) functorial cofibrant replacement in the category of simplicial spaces.
(At this point, some decisions must be made on the meaning of \emph{space} and on preferred model category
structures on the category of spaces and on the category of simplicial spaces. Let us say, for example, that
\emph{space} just means topological space and that we use the standard model category structure on the category
of spaces where \emph{fibration} means \emph{Serre fibration} and the weak equivalences are the maps which are
classically called weak equivalences. There is then a preferred model category structure on the category of simplicial
spaces where the weak equivalences are the simplicial maps which are degreewise weak equivalences in the category of spaces,
and the fibrations are degreewise fibrations in the category of spaces. This is a good choice for our purposes.)
\end{itemize}
It is well known that the classifying space $B\sC_X$ is then a correct model for
the spaces $B\haut_{N\fin}(X)$, where $\haut_{N\fin}(X)$ is the union of the homotopy invertible components of $\rmap_{N\fin}(X,X)$.
A diagram of the shape
\[  \xymatrix{
 X(0) & \ar[l]_-{v(1)}  X(1) & \ar[l]_-{v(2)}   \cdots & \ar[l]_-{v(k)}  X(k)
} \]
in $\sC_X$ determines a space $W(v(1),\dots,v(k))$, the inverse limit of the diagram of simplicial sets
\[ W(v(1))\to W(X(1)) \leftarrow W(v(2))\to W(X(2))
\leftarrow \cdots \to W(X(k-1)) \leftarrow W(v(k)). \]
This comes with forgetful projections to the $W(X(i))$ which are weak equivalences and
Kan fibrations. Now we have the following projection map.
\begin{equation} \label{eqn-reckless}
\hocolimsub{(v(1),\dots,v(k))} W(v(1),\dots,v(k))\qquad \lra  \hocolimsub{(v(1),\dots,v(k))}~\star
\end{equation}
These homotopy colimits are taken over the category where an object is a contravariant functor
from the ordered set $[k]=\{0,1,\dots,k\}$ to $\sC_X$~, for some $k$, and a morphism from $v\co [k]^\op\to \sC_X$
to $v'\co [\ell]^\op\to \sC_X$ is a monotone injective map $u\co [k]\to[\ell]$ such that $v'u=v$. The target of the
map~(\ref{eqn-reckless}) is
still an incarnation of $B\haut_{N\fin}(X)$ and the map itself is a quasi-fibration with fibers weakly
equivalent to $W(X)$. Therefore we can say that the map~(\ref{eqn-reckless}) determines a classifying map
\begin{equation} \label{eqn-hauthaut}  B\haut_{N\fin}(X) \lra B\haut(W(X))=B\haut(\rsec(\psi_X)) \end{equation}
where $\haut(W(X))$ is the union of the homotopy invertible path components of $\rmap(W(X),W(X))$.
\end{rem}

\begin{defn} \label{defn-shiftfunctorplus} Definition~\ref{defn-shiftfunctor} has a variant in which
$\finplus$ takes the place of $\fin$. This is straightforward. It begins with a
functor $\sigma\co \finplus\to \finplus$ given by disjoint sum with a
singleton. In more
detail, $\sigma$ is given by $[k]\mapsto [k+1]$ on objects, and for a morphism
$f\co [k]\to [\ell]$ (based map) the morphism $\sigma(f)$ is given by $\sigma(f)(x)=f(x)$ for $x\le k$
and $\sigma(f)(k+1)=\ell+1$. The standard inclusions of $[k]$ in $[k+1]$ define
a natural transformation $u\co \id\to \sigma$ between endofunctors of $\fin$\,. For a simplicial
space $X$ over $N\finplus$ we define $\psi_X\co \shift(X)\to X$ and $\rsec(\psi_X)$ as in
definition~\ref{defn-shiftfunctor}, mutatis mutandis.

The definition of an elementary Segal space over $N\fin$ carries over to $N\finplus$ almost without
change. We only have to describe $r$-simplices in $N\finplus$ as diagrams of \emph{based} finite sets
and based maps
\[   C = (\,[k_0] \xleftarrow{f_1} [k_1] \xleftarrow{f_2} [k_2]
\xleftarrow{f_3} \cdots \xleftarrow{f_r} [k_r]\,) \]
An element $x$ of $[k_j]$, where $j\in\{1,2,\dots,r-1\}$, is \emph{heavy } if its preimage under
$f_{j+1}$ has more than one element. No exception is made for $x=0$.
\end{defn}

\subsection{Back to occupants}
\label{subsec-Reedytriumph}
The homotopy limit of the contravariant functor $\Psi$ in \cite[\S2.1]{WeissSen1}
was defined using the Bousfield-Kan formula, i.e., as $\Tot$
of a certain cosimplicial space $[r]\mapsto \Gamma_r(\Psi)$. Here $\Gamma_r(\Psi)$ is the
section space of a fiber bundle $E_r^!(\Psi) \to N_r\sP(M\smin\partial M)$ such that the fiber over
\[  ((S_0,\rho_0)\ge (S_1,\rho_1)\ge \cdots\ge (S_r,\rho_r)) \]
is $M\smin V(S_r,\rho_r)$. \newline
 From a model category point of view we should have proceeded differently. The first step should have been to introduce
$E_r(\Psi)$, total space of a fiber bundle on $N_r\sP(M\smin\partial M)$ such that the fiber over
\[  ((S_0,\rho_0)\ge (S_1,\rho_1)\ge \cdots\ge (S_r,\rho_r)) \]
is $M\smin V(S_0,\rho_0)$. Note the difference between $E_r^!(\Psi)$ and $E_r(\Psi)$.
Now $[r]\mapsto E_r(\Psi)$ is a simplicial space and the projections
$E_r(\Psi)\to N_r\sP(M\smin\partial M)$ make up a simplicial map $p\co E(\Psi)\to N\sP(M\smin\partial M)$. It is easy to
see that $E(\Psi)$ is a complete Segal space like $N\sP(M\smin\partial M)$, although this will not be used
explicitly in the following. Think of $E(\Psi)$ as the
Grothendieck construction of the contravariant functor $\Psi$.
This suggests the definition
\begin{equation}  \label{eqn-grandholim}
\holim~\Psi := \rsec\big(p\co E(\Psi)\to N\sP(M\smin\partial M)\big). \end{equation}
Now we need to show that this is in agreement with the definition of $\holim~\Psi$ used in \cite[\S2.1]{WeissSen1}.
In this section we have favored the model structure on the category of simplicial spaces where a morphism
is a fibration if it is degreewise a Serre fibration (and a weak equivalence if it is degreewise a
weak equivalence). But the definition of $\holim~\Psi$ given in \cite[\S2.1]{WeissSen1} is more easily understood
in terms of the Reedy model structure on the category of simplicial spaces. It was already pointed out
in \cite[1.1.3]{WeissSen1} 
that the simplicial space $N\sP(M\smin\partial M)$ is Reedy cofibrant.
The map $p\co E(\Psi)\to N\sP(M\smin\partial M)$ is not (claimed to be) a Reedy fibration, but the
definition of $\holim~\Psi$ in \cite[\S2.1]{WeissSen1} contains a well-concealed suggestion
for a replacement by a Reedy fibration. Let
\[  x=((S_0,\rho_0)\ge (S_1,\rho_1)\ge \cdots\ge (S_r,\rho_r)) \]
be a point in $N_r\sP(M\smin\partial M)$. Let $F_{r,x}$ be the space of maps from
$\Delta^r$ to $M$ satisfying the condition that, for every monotone injective $u\co[t]\to[r]$,
the corresponding face of $\Delta^r$ is taken to $M\smin V(S_{u(t)},\rho_{u(t)})$. Let
\[  E^\sharp_r(\Psi) \to N_r(\sP(M\smin\partial M) \]
be the fibration such that the fiber over $x$ is $F_{r,x}$. There is an inclusion
\[ E_r(\Psi)\to E^\sharp_r(\Psi)~; \]
indeed $F_{r,x}\cap E_r(\Psi)$ is precisely the subspace of the
constant elements in $F_{r,x}$\,. Moreover it is easy to see that $E^\sharp_r(\Psi)$ is a simplicial space
again. In the factorization
\[  E(\Psi) \hookrightarrow E^\sharp(\Psi) \lra N\sP(M\smin\partial M) \]
of $p$\,, the first map is a weak equivalence and the second is a Reedy fibration. Therefore it is allowed
to define $\holim~\Psi$ as the space of sections of the map of simplicial spaces
\begin{equation} \label{eqn-Reedytriumph}
E^\sharp(\Psi) \lra N\sP(M\smin\partial M) \end{equation}
and this is exactly the definition of $\holim~\Psi$ given in \cite[\S2.1]{WeissSen1}. \newline
Now it is also clear how we can relate the older definition of $\holim~\Psi$
to the alternative definition~(\ref{eqn-grandholim}).
Namely, we pass from the honest section space of~(\ref{eqn-Reedytriumph}) to the derived section space
of~(\ref{eqn-Reedytriumph}) in a possibly different model category structure (with the same
weak equivalences), and compare that
to the derived section space~(\ref{eqn-grandholim}).

\begin{expl} \label{expl-goodthing} One thing that we must take away from this section is an identification (zigzag of weak equivalences) of
$\holim~\Psi$ in section~\ref{sec-intro} or \cite[\S2.1]{WeissSen1}
with $\rsec(\psi_X)$ of definition~\ref{defn-shiftfunctor}, where $X$ is $\config(M\smin\partial M)$.
There are a few simple steps to this conversion.
\begin{itemize}
\item[(i)] We start with $X=N\sP(M\smin\partial M)$ and the map $p\co E\to X$ of simplicial spaces
where $E_r=E_r(\Psi)$ is the total space of a bundle on $X_r$ such that the fiber over a point
$((S_0,\rho_0)\ge \cdots\ge (S_r,\rho_r))$ is $M\smin V(S_0,\rho_0)$.
By definition and by the foregoing discussion, $\holim~\Psi$ in section~\ref{sec-intro}
is $\rsec(p)$, which can be thought of as the space of sections of $p^\sharp\co E^\sharp\to X$,
where $p^\sharp$ is a Reedy fibration replacing $p$.
\item[(ii)] We modify $X$, $E$ and $E^\sharp$ in (i) by choosing total orderings for all configurations in sight.
The new $X$ is now entitled to the name $\config(M\smin\partial M)$ and it is a simplicial
space over $N\fin$. We obtain a new $\rsec(p)$, space of sections of the new $p^\sharp\co E^\sharp\to X$.
The space $\rsec(p)$ in (i) maps to the new version $\rsec(p)$ here in (ii) by a weak equivalence.
\item[(iii)] We keep $X$ as in (ii) but make some small changes to $E$ and $E^\sharp$. The new $E_r$ is the total space of a fiber
bundle on $X_r$ such that the fiber over a string $((S_0,\rho_0)\ge \cdots\ge (S_r,\rho_r))\in X_r$\,, where the sets $S_0,\dots,S_r$
are totally ordered, is the space of pairs $(z,\varepsilon)$ where $z\in M\smin V(S_0,\rho_0)$ and $\varepsilon$ is a positive
real number which is less
than the distance from $z$ to $\partial M$ and less than the distance from $z$ to the closure of $V(S_0,\rho_0)$.
The new $E$ is entitled to the name $\shift(X)$.
We obtain a new variant of $p\co E\to X$ which is entitled to the name $\psi_X\co \shift(X)\to X$.
We obtain a new  variant of $\rsec(p)$.
This is entitled to the name $\rsec(\psi_X)$ where $X=\config(M\smin\partial M)$.
\end{itemize}
Similar good things can be said about $\holim~\Psi$ in \cite[\S3.2]{WeissSen1}.
It can be identified with $\rsec(\psi_X)$ of definition~\ref{defn-shiftfunctorplus}, where $X$ is
$\config(M\smin\partial_1M)$, a fiberwise complete Segal space over $N\finplus$\,. To recall some
of the details: $M$ is a smooth compact Riemannian manifold with boundary and corners in the boundary,
so that
\[ \partial M=\partial_0M\cup\partial_1M \]
where $\partial\partial_0M=\partial\partial_1M=\partial_0M\cap\partial_1M$.
The topological poset $\sP(M\smin\partial_1M)$ has elements $(S,\rho)$ where $S$ is a finite subset of $M\smin\partial M$
and $\rho\co S\sqcup\partial_0M\to\RR$ is a function with positive values, locally constant on $\partial_0M$.
There are some smallness conditions on $\rho$ as usual. In addition it is required that the Riemannian metric on $M$ be
a product metric near $\partial_1M$ (product of a Riemannian metric on $\partial_1M$ and the standard metric on
a closed interval). For $(S,\rho)\in \sP(M\smin\partial M)$ the set $V(S,\rho)$ is defined as an open subset
of $M\smin\partial_1M$, but the functor $\Psi$ is defined by $(S,\rho)\mapsto M\smin V(S,\rho)$. As a result there is a
map
\[ \partial_1M\to \holim~\Psi \]
induced by the inclusion $\partial_1M\to M\smin V(S,\rho)$. This is a weak equivalence under some (rather severe) conditions
on $M$.
\end{expl}

\section{Occupants and homotopy automorphisms}  \label{sec-mainresult}
\subsection{Boundary of a compact manifold as a homotopy link}
For a smooth compact $M$, let $\holink(M/\partial M,\star)$ be the space of paths
$w\co [0,1]\to M/\partial M$ which satisfy $w^{-1}(\star)=\{0\}$, with the compact-open topology.
Evaluation at $1\in [0,1]$ gives a map
\begin{equation} \label{eqn-bdryholink}
q_M\co \holink(M/\partial M,\star) \lra M\smin\partial M\,.
\end{equation}
It is well known, and it will be made precise in a moment, that $q_M$ is a good homotopical substitute for the inclusion map
$\partial M\to M$.

\smallskip
The modest advantage that $q_M$ has for us, compared to the inclusion $\partial M\to M$, is that
the canonical action of the homeomorphism group $\homeo(M)$ on the map $q_M$ extends rather obviously to an action
of the homeomorphism group $\homeo(M\smin\partial M)$ on $q_M$
(because $M/\partial M$ is the one-point compactification of $M\smin\partial M$.)
Of course, $\homeo(M)$ also acts canonically on the inclusion $\partial M\to M$, or equivalently
on the pair $(M,\partial M)$, but this action does not extend obviously or otherwise to an action of
$\homeo(M\smin\partial M)$ except in the cases where $\partial M=\emptyset$ or $\dim(M)\le 1$.

We take the view that $\homeo(M\smin\partial M)$ is an enlargement of
$\homeo(M)$. Indeed, the restriction homomorphism from
$\homeo(M)$ to $\homeo(M\smin\partial M)$
is injective, due to the fact that $M\smin\partial M$ is dense in $M$.

\smallskip
Let $Z_M$ be the space of maps $w\co [0,1]\to M$ such that $w^{-1}(\partial M)=\{0\}$.
Then we have a map $Z_M\to \holink(M/\partial M,\star)$ given by composing
elements $w\in Z_M$ with the quotient map $M\to M/\partial M$. This map is a weak homotopy
equivalence. There is also a forgetful map $Z_M\to \partial M$ given by evaluation at $0$.
Together these maps make up a diagram
\[
\xymatrix@M=6pt@C=12pt{
\holink(M/\partial M,\star) \ar[rr]^-{q_M} && M\smin\partial M \ar[d]^-\simeq \\
Z_M \ar[u]^-\simeq \ar[r]^-\simeq & \partial M \ar[r] ^-{\textup{incl.}} & M
}
\]
which is commutative up to a preferred homotopy. The group $\homeo(M)$ acts on the
whole diagram, respecting the preferred homotopy. For the top row that action extends to
an action of $\homeo(M\smin\partial M)$.

\smallskip
To package some of this in a more memorizable way, let us write $M_-$ for
$M\smin\partial M$ and $\partial^h M_-$ for
$\holink(M/\partial M,\star)$. Then we have
\[  q_M\co \partial^h M_-\lra M_- \,. \]
As we have seen, this map is a good homotopical substitute for the inclusion of $\partial M$ in $M$.
The group $\homeo(M_-)$ acts on $q_M$.
The action determines a map
\[
B\homeo(M_-) \lra B\haut(q_M\co \partial^h M_-\to M_-)
\]
where $\haut(-)$ generally denotes a space of derived homotopy automorphisms of an object
in some model category. (A map such as $q_M$ should be viewed as a functor from the totally ordered set $\{0,1\}$
to spaces. We use one of the standard model category structures on the category of such functors.)

\subsection{The prototype} Suppose that the compact smooth $M$ satisfies the
condition of \cite[Thm.~2.1.1]{WeissSen1}. That is to say, $M$ is the total space of a smooth disk bundle
on a smooth compact manifold $L$ without boundary, where the fibers are of dimension $c\ge 3$. (The fibers are smooth manifolds homeomorphic to the disk $D^c$.) Write
\[  \haut_{N\fin}(\config(M\smin\partial M)) \]
for the union of the homotopy invertible (path) components in
$\rmap_{N\fin}(X,X)$ where $X=\config(M\smin\partial M)$. Composition makes this into
a grouplike topological or simplicial monoid.

\begin{thm} \label{thm-factabs} Under these conditions on $M$,
the broken arrow in the following homotopy commutative diagram can be supplied:
\[
\xymatrix@C=35pt@M=8pt@R=20pt{
B\homeo(M_-) \ar[r]^-{\textup{action}} & B\haut(\,\partial^h M\to M_-\,)  \\
B\homeo(M_-) \ar@{=}[u] \ar[r]^-{\textup{action}} & B\haut_{N\fin}(\config(M_-)) \ar@{..>}[u]
}
\]
\end{thm}

\proof Write $X$ for $\config(M_-)$; we use the particle model. Let $A=\config(M_-\,;0)$
be the simplicial subspace of $X$ obtained by allowing only configurations
of cardinality zero. This is of course rather trivial: $A_r$ is a point for all $r\ge 0$.
We will be interested in $\psi_X\co \shift(X)\to X$ and in the section space $\rsec(\psi_X)$ and also
in the section space $\rsec(\psi_X|_A)$. It has already
been indicated that $\rsec(\psi_X)$ is weakly equivalent to $\partial M$, and it is easy to show directly
that $\rsec(\psi_X|_A)$ is weakly equivalent to $M$ or to $M_-$\,. We need a more precise statement.
\begin{itemize}
\item[(i)] The action of $\homeo(M_-)$ on $q_M\co \partial^h M\to M_-$ gives rise to two
fiber bundles on $B\homeo(M_-)$, with fibers $\partial^h M$ and $M_-$ respectively,
and a map from the first to the second. To this we refer loosely as a \emph{pair} of fiber bundles
on $B\homeo(M_-)$.
\item[(ii)] The action of $\haut_{N\fin}(X)$ on $\rsec(\psi_X)\to \rsec(\psi_X|_A)$
gives rise to two fibrations on $B\haut_{N\fin}(X)$, with fibers $\rsec(\psi_X)$ and
$\rsec(\psi_X|_A)$ respectively, and a map from the first to the second.
To this we refer loosely as a \emph{pair} of fibrations on $B\haut_{N\fin}(X)$.
\item[(iii)] Under pullback along the inclusion
\[ B\homeo(M_-)\lra B\haut_{N\fin}(\config(M_-))=B\haut_{N\fin}(X),  \]
the fibration pair in (ii) becomes fiberwise homotopy equivalent to the fiber bundle pair in (i).
More precisely, there is a zig-zag of fiberwise weak equivalences over $B\homeo(M_-)$, etc.
\end{itemize}
What we have to prove is (iii). To keep notation manageable, we will concentrate
on the boundary fibrations with fibers $\rsec(\psi_X)$ and $\partial^h M$,
neglecting the fibrations with fibers $\rsec(\psi_X|_A)$ and  $M_-$\,.
The main ideas are already in place and we just arrange them
by introducing a simplicial map
\[  \varphi_X\co X^* \lra X~, \]
closely related to $\psi_X\co E^\sigma(X)\to X$. For $r\ge 0$, the space $X^*_r$
is an open subset of
\[ X_r\times\partial^h M_- \]
consisting of all pairs $(y,w)$ where $y\in X_r$ and $w\co [0,1]\to M/\partial M$ is an element of
$\partial^hM_-$ such that $w(t)$ for $t>0$ is not contained in the \emph{support} of
$y$. (The support of $y$ is a compact subset of $M_-$. It is the union of all finite subsets which arise as images
of the various maps from finite sets to $M_-$ which make an appearance in  the description of $y$.)
\begin{itemize}
\item[(iv)] The inclusion $X^* \to X\times \partial^h M_-$ is a degreewise weak equivalence
(where $\partial^h M_-$ should be viewed as a constant simplicial space).
\item[(v)] There is a map $X^*\to E^\sigma(X)$ over $X$, given by taking a pair
$(y,w)\in X_r$ to the element of $E^\sigma_r(X)$ obtained by evaluating $w$ at $1$ and using that value
to increase the cardinality of all configurations in $y$ by one.
\item[(vi)] That map $X^*\to E^\sigma(X)$ over $X$ induces a weak equivalence of derived
section spaces,
\[  \rsec(\varphi_X\co X^*\to X) \lra \rsec(\psi_X\co E^\sigma(X)\to X) \]
by \cite[Thm.~2.1.1]{WeissSen1}, translated to the configuration category setting.
\item[(vii)] Projection from $X^*$ to $\partial^h M_-$ determines a weak equivalence
\[  \rsec(\varphi_X\co X^*\to X)\lra \partial^h M_-\,.  \]
\end{itemize}
The weak equivalences in (vi) and (vii) respect canonical actions of $\homeo(M_-)$.
Equivalently, they extend to fiberwise weak homotopy equivalences between fibrations over $B\homeo(M_-)$.
This gives us the zig-zag of fiberwise weak homotopy equivalences which we require. \qed

\begin{rem} S.~Tillmann in \cite[Thm 1.2]{TillSimp} makes a statement stronger than \cite[Thm.~2.1.1]{WeissSen1}.
The conclusion is the same but the condition on $M$
is weaker. Namely, the manifold $M$ is required to be a \emph{smooth thickening} of an embedded
simplicial complex $K$ of codimension at least $3$. This appears to be the right level of generality.
See \cite{TillSimp} for the precise condition on $M$. Using this stronger form,
one can deduce a stronger form of theorem~\ref{thm-factabs} with the corresponding weaker condition on $M$.
(At the time of writing, the status of \cite{TillSimp} is \emph{submitted}.)
\end{rem}

\subsection{Variant with cardinality restriction} \label{subsec-factabstrunc}
For a more technical variant of theorem~\ref{thm-factabs} we need the Postnikov decomposition of the map
$\partial M\to M$. For an integer $a\ge 0$ there is the factorization
\[  \partial M \lra \wp_a \partial M \lra M \]
of $\partial M\to M$ where $\wp_a\partial M$ is obtained from $\partial M$, as a space over $M$,
by killing the relative homotopy groups of $\partial M\to M$ in dimensions $\ge a+2$. So
$\partial M\to \wp_a\partial M$ is $(a+1)$-connected. We apply this construction (and use
the informal notation) mutatis mutandis with $\partial^h M_-\to M_-$ instead of $\partial M\to M$.

For $M$ satisfying the condition of \cite[Thm.~2.1.1]{WeissSen1} and an integer $j\ge 1$, let
$\config(M_-\,;j)$ be the truncated configuration category where only configurations
of cardinality $\le j$ are allowed. The condition on $M$ is, imprecisely stated, that $M$ is the total space
of a smooth disk bundle of fiber dimension $c\ge 3$ on a smooth closed manifold.
The truncated configuration category $\config(M_-\,;j)$ is still a fiberwise complete Segal space
over $N\fin$.

\begin{thm} \label{thm-factabstrunc} Under these conditions on $M$,
the broken arrow in the following homotopy commutative diagram can be supplied:
\[
\xymatrix@C=35pt@M=8pt@R=20pt{
B\homeo(M_-) \ar[r]^-{\textup{action}} & B\haut(\,\wp_{(j+1)(c-2)}\partial^h M_-\lra M_-\,)  \\
B\homeo(M_-) \ar@{=}[u] \ar[r]^-{\textup{action}} & B\haut_{N\fin}\big(\config(M_-\,;j)\big) \ar@{..>}[u]
}
\]

\end{thm}

The proof of this follows the lines of the proof of theorem~\ref{thm-factabs} but relies more on the
connectivity estimates in \cite[Thm.~2.1.1]{WeissSen1}.  \qed

\section{Occupants and homotopy automorphisms, relative case}
As in \cite[\S4.1]{WeissSen1}, let  
$M$ be a smooth compact manifold with boundary and corners, so that $\partial M$ is the union of
smooth codimension zero submanifolds $\partial_1M$ and $\partial_0M$ satisfying
$\partial\partial_0M=\partial\partial_1M =\partial_0M\cap \partial_1M$.
The group of homeomorphisms of $M\smin\partial_1M$
which restrict to the identity on $\partial_0M\smin\partial_1M$ acts in an $A_\infty$ sense, and by homotopy automorphisms,
on the square of inclusion maps
\[
\xymatrix@R=12pt{
\partial_0M\cap \partial_1M \ar[d] \ar[r] & \ar[d] \partial_0M \\
\partial_1 M \ar[r] & M
}
\]
fixing the spaces in the top row (i.e., the pair $(\partial_0M,\partial_0M\cap \partial_1M)$) pointwise.
We look for an extension of this
action to (or factorization of this action through)
\[ \haut_{N\finplus}(\config(M\smin\partial_1 M);\config(U\smin\partial_1M)) \]
where $U$ is a collar neighborhood of $\partial_0M$ in $M$. The semicolon notation means that we
allow only derived automorphisms of $\config(M\smin\partial_1 M)$ which extend the identity on
$\config(U\smin\partial_1M)$. We find such an extension or factorization under severe conditions. The proofs are similar to those in
sections~\ref{sec-shifty} and~\ref{sec-mainresult}, but technically more
demanding. We need the twisted arrow construction on $\config(M\smin\partial_1M)$ and we need to reformulate the
construction in \cite[\S4.1]{WeissSen1} of a functor $\Theta$ from there to spaces in a more general setting.

\subsection{Twisted arrow construction on simplicial spaces}
The twisted arrow construction on a simplicial space $X$ is $\tw(X):=X\circ\beta$,
where $\beta\co \Delta\to \Delta$ is the functor $[n]\mapsto[2n+1]$. More precisely, $\Delta$ is the category of totally ordered nonempty finite sets
and order-preserving maps, or the equivalent full subcategory with objects $[n]$ for $n\ge 0$, and
$\beta$ is the functor which takes a totally ordered
nonempty finite set $S$ to $S\sqcup S^\op$ (with the concatenated total ordering where $a<b$ if $a\in S\subset S\sqcup S^\op$
and $b\in S^\op\subset S\sqcup S^\op$). \newline
The inclusions $S\to S\sqcup S^\op$ define a natural transformation $e\co \id\to \beta$. This induces a simplicial map
$\tw(X)\to X$ of simplicial spaces which is entitled to the name \emph{source}. Example: If $X=N\sC$ for a small category $\sC$, then $\tw(X)=N(\tw(\sC))$ where $\tw(\sC)$
is the twisted arrow category of $\sC$. (An object of $\tw(\sC)$ is a morphism in $\sC$ and a morphism in $\tw(\sC)$ is a commutative diagram
\[
\xymatrix@R=10pt{ a \ar[d] \ar[r] & b \\
c \ar[r] & d \ar[u]
}
\]
in $\sC$, where the top row is the source object in $\tw(\sC)$ and the bottom row is the target object.) The canonical
map $\tw(X)\to X$ is then the map of nerves induced by the forgetful functor which takes an object in $\tw(\sC)$, alias
morphism $a\to b$ in $\sC$, to its source $a$.

\subsection{More shifting} \label{subsec-twistshift} As mentioned in definition~\ref{defn-shiftfunctorplus},
the functor $\sigma$ from $\fin$ to $\fin$ and the natural transformation $u\co \id\to \sigma$
of definition~\ref{defn-shiftfunctor}
extend in a straightforward way to a functor $\finplus \to \finplus$
and a natural transformation from $\id\co \finplus\to \finplus$ to $\sigma$. These are still denoted $\sigma$ and $u$,
respectively. Here we need another variant consisting of a functor $\tau\co \tw(\finplus)\to \tw(\finplus)$
and a natural transformation $v$ from the identity on $\tw(\finplus)$ to $\tau$.

\begin{defn}  On objects,
$\tau\co \tw(\finplus)\to \tw(\finplus)$ is defined by
\[  \tau\big(f\co [m] \to [n]\big):= \big(g\co [m+1] \to [n]\big) \]
where $g(x)=f(x)$ for $x\in [m]$ and $g(m+1)=0$. The remaining details are settled in such a way that
$F_s\tau=\sigma F_s$ and $F_t\tau=F_t$~, where $F_s,F_t\co \tw(\finplus)\to \finplus$ are the
functors given by source and target, respectively. \newline
For a fiberwise complete Segal space $Y$ over $\tw(N\finplus)=N(\tw(\finplus))$ let $\shiftt(Y)$ be the Segal
space defined by the pullback square of simplicial spaces
\[
\xymatrix{  \shiftt(Y) \ar[d] \ar[r] &  Y^{\Delta[1]} \ar[d] \\
N\tw(\finplus) \ar[r]^-{(\tau,v)} &  (N\tw(\finplus))^{\Delta[1]}_{\rule{0mm}{1mm}}
}
\]
There is a map $\theta_Y\co \shiftt(Y)\to Y$ over $N\tw(\finplus)$ given by composing the
map $\shiftt(Y)\to Y^{\Delta[1]}$
from the defining pullback square with the map
\[ Y^{\Delta[1]}\to Y^{\Delta[0]}\cong Y \]
determined by the
map $\Delta[0]\to \Delta[1]$ which takes the preferred generator in degree $0$ to $d_0$
of the preferred generator in degree $1$.
\end{defn}

\begin{prop} Let $X$ be a Segal space over $N\finplus$ and let $Y=\tw(X)$, Segal space over
$\tw(N\finplus)=N\tw(\finplus)$. If $X$ is elementary over $N\finplus$, then $\theta_Y$ is a
right fibration.
\end{prop}

\proof This follows the lines of the proof of proposition~\ref{prop-elementaryright}.
Let $f=(f_0,f_1)$ be a morphism in $\tw(\finplus)$ given by a diagram
\[
\xymatrix{
[k_0]  & \ar[l]_-{g_0} [k_1] \ar[d]^-{f_1} \\
[\ell_0] \ar[u]^-{f_0} & \ar[l]^-{g_1}  [\ell_1]
}
\]
The main point is to show
that the face operator $d_1\co Y_2\to Y_1$ restricts to a weak equivalence $Y_2(A'')\to Y_1(d_1A'')$,
where $A''$ is a certain $2$-simplex in $N\tw(\finplus)$ determined by $f$; this takes the place of $A''$ in the proof
of proposition~\ref{prop-elementaryright}. Here the diagram $A''$ is
\[
\begin{aligned}
\xymatrix{
{\big([k_0]\xleftarrow{g_0}[k_1]\big)}  \ar@{..>}[dr]  \ar[d]_-{v}  \\       
{\tau\big([k_0]\xleftarrow{g_0}[k_1]\big)} \ar[r]_-{\tau(f)}  & {\tau\big([\ell_0]\xleftarrow{g_1}[\ell_1]\big)} \\
}
\end{aligned}
\]
By the definition of $Y$ as $\tw(X)$, the map $Y_2(A'')\to Y_1(d_1A'')$ induced by $d_1$ can be
identified with the map $X(C)\to X(d_1d_4 C)$ induced by $d_1d_4\co X_5\to X_3$, where $C$ is certain
$5$-simplex in $N\finplus$. In this way the elementary property of $X$ can be exploited. \qed

\bigskip
Let $Y$ be a fiberwise complete Segal space over $N\tw(\finplus)$.
Reasoning as in section~\ref{subsec-shifty}, we find that the space of derived homotopy automorphisms
$\haut_{N\tw(\finplus)}(Y)$ (which is the union of some path components of $\rmap_{N\tw(\finplus)}(Y,Y)$)
acts on the derived section space $\rsec(\theta_Y)$. This action is constructed as a map
\[  B\haut_{N\tw(\finplus)}(Y) \lra B\haut(\rsec(\theta_Y)). \]
More specifically, if $Y=\tw(X)$ for a fiberwise complete Segal space $X$ over $N\finplus$, then we also have an
obvious map $B\haut_{N\finplus}(X)\to B\haut_{N\tw(\finplus)}(Y)$ which we can compose with the above
to conclude that $\haut_{N\finplus}(X)$ acts on $\rsec(\theta_Y)$.

\begin{expl} \label{expl-goodthing2}
Of particular interest is the case where $Y=\tw(X)$ and $X$ is the simplicial space $\config(M\smin\partial_1M)$, for
a smooth compact manifold $M$ with corners, $\partial M=\partial_0M\cup\partial_1M$ etc.,
as in \cite[\S4.1]{WeissSen1}. In that case the derived section
space $\rsec(\theta_{\tw(X)})$
can be identified with the space $\holim~\Theta$ in \cite[\S4.1]{WeissSen1}. The reasoning
is analogous to that in example~\ref{expl-goodthing}. The commutative diagram
\begin{equation} \label{eqn-oldloc}
\begin{aligned}
\xymatrix@C=16pt{
\partial\partial_1M \ar[rr]^-{\textup{inclusion}} \ar[d] && \partial_1M \ar[d] \\
\holim~\Theta \ar[r]  & \holim~\Psi\circ F_s & \ar[l]_-\simeq  \holim~\Psi
}
\end{aligned}
\end{equation}
of \cite[\S4.1]{WeissSen1} can be recast as
\begin{equation} \label{eqn-oldloc2}
\begin{aligned}
\xymatrix@C=16pt{
\partial\partial_1M \ar[rr]^-{\textup{inclusion}} \ar[d] && \partial_1M \ar[d] \\
\rsec(\theta_{\tw(X)}) \ar[r]  & \rsec(F_s^*\psi_X) & \ar[l]_-\simeq \rsec(\psi_X)
}
\end{aligned}
\end{equation} 
where $F_s^*\psi_X$ is defined by a (degreewise) homotopy pullback square of simplicial spaces
and simplicial maps
\[
\xymatrix{
F_s^*\shift(X) \ar[d]_-{F_s^*\psi_X} \ar[r]  &  \shift(X) \ar[d]^-{\psi_X} \\
   \tw(X) \ar[r]^-{F_s} & X\,.
}
\]
It is a key point that $\holim~\Theta$ has a
locality property. Briefly, its homotopy type depends only on an arbitrarily small neighborhood $U$ of $\partial_0M$
in $M$. This is established in section~\ref{sec-local} by reduction to a weaker statement of that type
proved in \cite[\S4.2]{WeissSen1}. Doubts about the stronger form raised in \cite[\S4.2]{WeissSen1}
have turned out to be unjustified.
\end{expl}

\subsection{Variant of the main result with gate} \label{subsec-gatefact}
In theorems~\ref{thm-factrel} and~\ref{thm-factreltrunc} below
some unsystematic notation is used to describe spaces of automorphisms of pairs and more complicated
situations. Fix $M$, a smooth manifold with boundary and corners, so that
$\partial M=\partial_0M\cup \partial_1M$
and $\partial\partial_0M=\partial\partial_1M=\partial_0M\cap\partial_1M$. We write,
rather unsystematically,
\[  M_-:= M\smin \partial_1M\,, \qquad \partial_0 M_-= \partial_0M\smin\partial\partial_0M \]
and $\partial^h_1M_-$ for $\holink(M/\partial_1M,\star)$, as well as
$\partial^h\partial_0M_-$ for $\holink(\partial_0M/\partial\partial_0M,\star)$. There is
a commutative square
\[
\xymatrix{ \partial^h\partial_0M_- \ar[r] \ar[d]^-{\textup{incl.}} & \partial_0M_- \ar[d]^-{\textup{incl.}} \\
\partial^h_1M_- \ar[r] & M_-
}
\]
where the horizontal maps are given by evaluation of paths at time $1$. This is our
preferred homotopical substitute for the square of inclusion maps
\[
\xymatrix{ \partial\partial_0M \ar[r] \ar[d] & \partial_0M \ar[d] \\
\partial_1M \ar[r] & M
}
\]
and the advantage of the substitute over the original is that the homeomorphism group $\homeo(M_-;\partial_0M_-)$
acts in a canonical way. (More unsystematic notation here: $\homeo(M_-;\partial_0M_-)$ consists of homeomorphisms
$M_-\to M_-$ which restrict to the identity on the boundary $\partial_0M_-$\,.)

One more abbreviation: let $U$ be a standard open collar neighborhood of $\partial_0M$ in $M$~. We assume
that the closure of $U$ in $M$ is a smooth closed collar. 
Write $U_-:= U\cap M_-$\,. Although it is a little careless, where
homeomorphisms $M_-\to M_-$ are mentioned which restrict to the identity on $\partial_0M_-$~, we may mean
homeomorphisms $M_-\to M_-$ which restrict to the identity on all of $U_-$.

\begin{thm} \label{thm-factrel} If $M$ satisfies the conditions of \cite[Thm.~3.2.1]{WeissSen1},
then the broken arrow in the following homotopy commutative diagram can be supplied:
\[
\xymatrix@C=40pt@M=6pt{
B\homeo(M_-;\,\partial_0M_-) \ar@{=}[d] \ar[r]^-{\textup{action}}
& {B\haut\begin{pmatrix}                                           
\partial^h\partial_0M_- & \!\!\!\lra\!\!\!&  \partial_0M_- & &  \\
\downarrow & &  \downarrow & \!\! ; \!\!& \textup{top row}  \\
\partial^h_1M_- & \!\!\!\lra\!\!\! &  M_- & &
\end{pmatrix}} \\
B\homeo(M_-;\,\partial_0M_-) \ar[r]^-{\textup{action}} &
B\haut_{N\finplus}(\config(M_-)~;~\config(U_-)) \ar@{..>}[u]
}
\]
\end{thm}

For clarification, this theorem has theorem~\ref{thm-factabs} as a special case. It is the special case where
$\partial_0M$ is empty.

\proof[Outline of a proof] Although the proof is similar to the proof
of theorem~\ref{thm-factabs}, it does require and use one additional idea. The outline will concentrate
on that.

Let $X=\config(M_-)$ and $\partial X:= \config(U_-)$, both
to be viewed as a complete Segal spaces over $N\finplus$\,.
As indicated earlier we can pretend that $\homeo(M_-;\partial_0M_-)$ consists of the
homeomorphisms $M_-\to M_-$ which restrict to the identity on (the closure of) $U_-$\,.
This is necessary to make the lower horizontal arrow
in the square (of the theorem) meaningful. \newline
We proceed initially as in the proof of theorem~\ref{thm-factabs}. In particular we have $A\subset X$ as before.
We are guided by the idea that the restriction map
\[  \rsec(\psi_X) \to \rsec(\psi_X|_A) \]
is weakly equivalent to the inclusion map
$\partial_1M\to M$, or to the homotopical substitute
\[ \partial^h_1M_- \lra M_-\,. \]
(This follows easily from \cite[Thm.~3.2.1]{WeissSen1} and the partial reformulation
at the very end of example~\ref{expl-goodthing}.) But we are chiefly interested in
homotopy automorphisms of the square
\begin{equation} \label{eqn-prehybrid}
\begin{aligned}
\xymatrix{
\partial^h\partial_0M_- \ar[r]\ar[d] & \partial_0M_- \ar[d] \\
\partial^h_1M \ar[r] & M_-
}
\end{aligned}
\end{equation}
fixing the top row pointwise. Therefore it seems that
we need to come to terms with a semi-combinatorial analogue of diagram~(\ref{eqn-prehybrid})
in the shape of a square
\begin{equation} \label{eqn-hybrid}
\begin{aligned}
\xymatrix@R=15pt{
\partial^h\partial_0M_- \ar[r]\ar[d] & \partial_0M_- \ar[d] \\
\rsec(\psi_X) \ar[r] & \rsec(\psi_X|_A).
}
\end{aligned}
\end{equation}
This is in agreement with the proof of theorem~\ref{thm-factabs}; there we had
$\partial_0M_-=\emptyset$. In general, we do \emph{not} (yet) have a sufficiently well understood combinatorial
expression for $\partial^h\partial_0M_-$ in terms of the configuration categories
$X$ and/or $\partial X$. Therefore diagram~(\ref{eqn-hybrid}) is the hybrid that it is. ---
The new and perhaps slightly unexpected task therefore is to set up diagram~(\ref{eqn-hybrid}) in such a way
that the actions of $\aut_{N\finplus}(X;\partial X)$ and $\haut_{N\finplus}(X;\partial X)$ on the lower row extend
to actions on the entire square which are trivial on the terms in the top row.
One solution is to construct diagram~(\ref{eqn-hybrid}) as the contraction of a bigger commutative diagram
\begin{equation} \label{eqn-hybrid2}
\begin{aligned}
\xymatrix@R=19pt{
\partial^h\partial_0M_- \ar[r]\ar[d]^a & \partial_0M_- \ar[d] \\
{\rsec(\theta_{\tw(X)})} \ar[r]\ar@{..>}[d] & {\rsec(\theta_{\tw(X)}|_{\tw(A)})} \ar@{..>}[d] \\
{\rsec(\psi_X)} \ar[r] & {\rsec(\psi_X|_A)}.
}
\end{aligned}
\end{equation}
where we use the notation of section~\ref{subsec-twistshift}.
Diagram~(\ref{eqn-oldloc2}) provides the dotted vertical arrows.
The lower square in diagram~(\ref{eqn-hybrid2}) is combinatorial, i.e., expressed in homotopical terms of $X$ and $\partial X$.
Then $\haut_{N\finplus}(X;\partial X)$
can act on the lower square. The action on the upper row (of the lower square) can be trivialized; this is made possible by
\cite[Prop.~4.2.1]{WeissSen1} and lemma~\ref{lem-longawait} below. Those actions can then
be extended canonically to actions on the entire diagram~(\ref{eqn-hybrid2})
which are trivial on the terms in the top square. ---
From this point onwards, the proof is a straightforward adaptation of the proof of theorem~\ref{thm-factabs}. \qed

\smallskip
The arrow labeled $a$ in diagram~(\ref{eqn-hybrid2}) is not claimed to be a weak equivalence. There is a suggestion
in \cite[\S4.2]{WeissSen1} that it is a weak equivalence under some additional
geometric hypotheses on $M$.

\subsection{Variant with gate and cardinality restriction}
We keep the notation of section~\ref{subsec-gatefact} and combine with the notation of section~\ref{subsec-factabstrunc}
for truncated configuration categories and Postnikov decompositions. Specifically (and unsystematically),
\[   \wp_a\partial^h_1M_- \lra  M_- \]
is the map obtained from $\partial^h_1M_- \lra M_-$ by killing the relative homotopy groups
in dimensions $\ge a+2$. There is still a commutative square
\[
\xymatrix{ \partial^h\partial_0M_- \ar[r] \ar[d] & \partial_0M_- \ar[d]  \\
\wp_{(j+1)(c-2)}\partial^h_1M_- \ar[r] &  M_-
}
\]

\begin{thm} \label{thm-factreltrunc} If $M$ satisfies the conditions of \cite[Thm.~3.2.1]{WeissSen1},
then the broken arrow in the following homotopy commutative diagram can be supplied:
\[
\xymatrix@C=33pt@M=5pt{
B\homeo(M_-;\,\partial_0M_-) \ar@{=}[d] \ar[r]^-{\textup{action}}
& {B\haut\begin{pmatrix}
\partial^h\partial_0M_- & \!\!\!\lra\!\!\!&  \partial_0M_- & &  \\
\downarrow & &  \downarrow \!\!\!&\!\! ; & \!\!\!\textup{top row}  \\
\wp_{(j+1)(c-2)}\partial^h_1M_- & \!\!\!\lra\!\!\! &  M_- & &
\end{pmatrix}} \\
B\homeo(M_-;\,\partial_0M_-) \ar[r]^-{\textup{action}} &
B\haut_{N\finplus}(\config(M_-;j)~;~\config(U_-;j)) \ar@{..>}[u]
}
\]
\end{thm}

This theorem has theorem~\ref{thm-factabstrunc} as a special case,
the case where $\partial_0M=\emptyset$
and consequently $\partial_1M=\partial M$.

\begin{appendices}
\section{Locality property of $\holim~\Theta$} \label{sec-local}  
\setcounter{thm}{0}
\renewcommand{\thesubsection}{\thesection}
We use the notation of \cite[\S4]{WeissSen1}.
So $M$ is a smooth compact manifold
with boundary and corners in the boundary, $\partial M=\partial_0M\cup\partial_1M$ etc., and
$\sP$ is short for $\sP(M\smin\partial_1M)$. For an object $(S,\rho)$ of $\sP$
let $V_\col(S,\rho)$ be the collar part of the open subset $V(S,\rho)$ in $M\smin\partial_1M$.
The contravariant functor $\Theta$ is defined on $\tw(P)$ by
\[  \Theta((S,\rho)\le (T,\sigma)):=\big(\textup{closure in $M$ of $V_\col(S,\rho)$}\big)\smin V(S,\rho). \]
A topological poset $\sQ$ is defined \cite[\S4.2]{WeissSen1} as a quotient
of $\tw(\sP)$. Two elements of $\tw(\sP)$, say
$\big((S,\rho)\le (T,\sigma)\big)$ and $\big((S',\rho')\le (T',\sigma')\big)$,
determine the same element of $\sQ$ if and only if $V_\col(T,\sigma)=V_\col(T',\sigma')$
and
\[ V_\col(T,\sigma)\cap V(S,\rho) =  V_\col(T',\sigma')\cap V(S',\rho'). \]
Therefore every element of $\sQ$ has a unique representative in $\tw(\sP)$ of the form
$((S,\rho)\le (\emptyset,\sigma))$.
For elements of $\sQ$ represented in this way by $((S,\rho)\le (\emptyset,\sigma))$
and $((S',\rho')\le (\emptyset,\sigma'))$, respectively, the first is $\le$
the second in $\sQ$ if and only if $V(\emptyset,\sigma')\subset V(\emptyset,\sigma)$
and $V(S,\rho)\cap V(\emptyset,\sigma')\subset V(S',\rho')$.

The name of the quotient functor from $\tw(\sP)$ to $\sQ$ is $K$. It is clear that
$\Theta=\Theta_1\circ K$ for a unique $\Theta_1$ from $\sQ$ to spaces.

\begin{lem} \label{lem-longawait} The map $\holim~\Theta_1\lra \holim~\Theta_1\circ K = \holim~\Theta$
induced by the forgetful functor $K\co \tw(\sP)\to \sQ$
is a weak equivalence.
\end{lem}

\proof There is a routine reduction (which is skipped here)
to the discrete setting; see \cite[\S1.2]{WeissSen1}. Let $\delta\sP$ and $\delta\sQ$ be the discrete posets obtained
from $\sP$ and $\sQ$,
respectively. Now we need to show that the canonical map
\[   \holim~\Theta_1|_{\delta\sQ}\lra \holim~\Theta_1K|_{\tw(\delta\sP)}  \]
is a weak equivalence. We start by introducing a poset $\sU$ intermediate between $\tw(\delta\sP)$ and $\delta\sQ$.
This is also a quotient of $\tw(\delta\sP)$. Two elements of $\tw(\delta\sP)$, say
$((S,\rho)\le (T,\sigma))$ and $((S',\rho')\le (T',\sigma'))$,
determine the same element of $\sU$ if and only if $V_\col(T,\sigma)=V_\col(T',\sigma')$
and $(S,\rho)=(S',\rho')$. For elements of $\sU$
represented in this way by $((S,\rho)\le (T,\sigma))$ and
$((S',\rho')\le (T',\sigma'))$ in $\tw(\delta\sP)$, the first is $\le$ the second if and only if
\begin{itemize}
\item $V(S,\rho)\subset V(S',\rho')$;
\item $V_\col(T',\sigma')\subset V_\col(T,\sigma)$;
\item $V(S',\rho')$ is in \emph{general position} to $V_\col(T,\sigma)$. This means that
for every connected component $W$ of $V(S',\rho')$, either $W\subset V_\col(T,\sigma)$ or
the closure of $W$ has empty intersection with the closure of $V_\col(T,\sigma)$.
\end{itemize}
Now the forgetful functor $K$ (in the discrete setting) is the
composition of two forgetful functors $K_s$ and $K_t$\,:
\[
\xymatrix{
\tw(\delta\sP) \ar[r]^-{K_t} &  \sU \ar[r]^-{K_s} & \delta\sQ\,.
}
\]
The plan is to show that both induced maps
\[  \holim~\Theta_1^\delta \lra \holim~\Theta_1^\delta K_s \lra \holim~\Theta_1^\delta K_sK_t \]
are weak equivalences (where $\Theta_1^\delta$ is $\Theta_1$ restricted to $\delta\sQ$).
It suffices to establish the property \emph{homotopy terminal} 
for the forgetful functors $K_t\co \tw(\delta\sP)\to \sU$
and $K_s\co \sU\to\delta\sQ$.   

\smallskip
\emph{Showing that $K_t$ is homotopy terminal.} Fix an object $z$ of $\sU$.
We can represent this by an object of $\tw(\delta\sP)$, say
\[  ((S^\sharp,\rho^\sharp)\le (T^\sharp,\sigma^\sharp)). \]
The category $(z\downarrow K_t)$
is identified with a full sub-poset $\sA$ of $\tw(\delta\sP)$, consisting of all objects
\[ ((S,\rho)\le (T,\sigma)) \in \tw(\delta\sP) \]
such that $(S^\sharp,\rho^\sharp)\le (S,\rho)$ in $\delta\sP$ and $V_\col(T,\sigma)\subset V(T^\sharp,\sigma^\sharp)$,
and $V(S,\rho)$ is in general position to $V_\col(T^\sharp,\sigma^\sharp)$. (See the earlier description of $\sU$.)
The poset $\sA$ has a full sub-poset $\sB$ consisting of all
$((S,\rho)\le (T,\sigma)) \in \sA$ where
\[ (S,\rho)=(S^\sharp,\rho^\sharp). \]
For every $y\in \sA$, the full sub-poset of $\sB$ consisting of the elements of $\sB$ which
are $\le y$ in $\sA$ has
a unique maximum. Indeed, if $y=((S,\rho)\le (T,\sigma))$ as above then that maximum
is $((S^\sharp,\rho^\sharp)\le (T,\sigma))$. Equivalently, the inclusion $\sB\to \sA$ has a right adjoint.
Now $\sB$ has a maximal element
given by $((S^\sharp,\rho^\sharp)\le (S^\sharp,\rho^\sharp))$. Therefore the
classifying space of $\sB$ is contractible, and so the classifying space of $\sA$ is also contractible.

\smallskip
\emph{Showing that $K_s$ is homotopy terminal.}
Fix an object $x$ of $\delta\sQ$ represented by
$((S^\flat,\rho^\flat)\le (T^\flat,\sigma^\flat))$ in $\tw(\delta\sP)$ in such a way that $T^\flat=\emptyset$.
We need to show that the poset $(x\downarrow K_s)$ has a
contractible classifying space. Identify $(x\downarrow K_s)$ with the full sub-poset $\sD$ of $\sU$ consisting of
all objects of $\sU$ represented by
\[ \big((S,\rho)\le (T,\sigma)\big)\in \tw(\delta\sP) \]
which 
satisfy
\begin{itemize}
\item[(i)] $V_\col(S^\flat,\rho^\flat)\subset V_\col(T,\sigma)\subset V(T^\flat,\sigma^\flat)$;
\item[(ii)] $V(S,\rho)\cap V_\col(T,\sigma)~\supset~V(S^\flat,\rho^\flat)\cap V_\col(T,\sigma)$.
\end{itemize}
Let $\sE$ be the full sub-poset of $\sD$ consisting of all objects as above which instead of (ii) satisfy the
stronger condition
\begin{itemize}
\item[(iii)] $V(S,\rho)\cap V_\col(T,\sigma)=V(S^\flat,\rho^\flat)\cap V_\col(T,\sigma)$.
\end{itemize}
For every $y\in \sD$ the full poset of the elements of $\sE$ which are $\le y$ in $\sD$ has
a unique maximum. Equivalently, the inclusion $\sE\to \sD$ has a right adjoint. It remains to show
that $\sE$ has a contractible classifying space. We show this in a separate step.

\smallskip
\emph{Showing that $\sE$ has a contractible classifying space.}
Let $\mathcal J$ be the full sub-poset of $\delta\sP$ consisting of the $(T,\sigma)\in\delta\sP$ which
have $T=\emptyset$ and satisfy the following additional conditions:
\begin{itemize}
\item[-] $V_\col(S^\flat,\rho^\flat)\subset V(T,\sigma)\subset V_\col(T^\flat,\sigma^\flat)$;
\item[-] $V(T,\sigma)=V_\col(T,\sigma)$ is in general position to $V(S^\flat,\rho^\flat)$.
\end{itemize}
Let $G$ be the functor from $\mathcal J^\op$ to posets which
\begin{itemize}
\item[-] to an object $(T,\sigma)\in \mathcal J$ associates
the full sub-poset of $\delta\sP$ consisting of all $(S,\rho)\in\delta\sP$ such that $V_\col(S,\rho)=V(T,\sigma)$;
\item[-] to a morphism $(T,\sigma)\le (T',\sigma')$ in $\mathcal J$ associates the map of posets
\[ G(T',\sigma') \ni (S',\rho') \mapsto (S,\rho) \in G(T,\sigma) \]
where $(S,\rho)$ is defined in such a way that $V_\col(S,\rho)=V(T,\sigma)$ and
\[ V(S,\rho)\smin V_\col(S,\rho) \]
is the (disjoint) union of $V(S',\rho')\smin V_\col(S',\rho')$ and the part of $V(S^\flat,\rho^\flat)$
contained in $V(T',\sigma')\smin V(T,\sigma)$.
\end{itemize}
Each $G(T,\sigma)$ has a contractible classifying space; indeed it has a minimal element.
It is also easy to see that
$\mathcal J^\op$ has a contractible classifying space. Therefore the
classifying space of the Grothendieck construction $\smallint G$
is contractible, e.g. by the Thomason homotopy colimit theorem \cite{Thomason}. But $\smallint G$ is clearly equivalent
to $\sE$.
(For the purposes here, by the Grothen\-dieck
construction $\smallint F$ of a functor $F$ from a small category $\mathcal M$ to
the category of small categories we mean the following category. Objects are pairs $(m,v)$
where $m$ is an object of $\mathcal M$ and $v$ is an object of $F(m)$. A morphism
from $(m,v)$ to $(n,w)$ is a pair $(f,g)$ where $f\co m\to n$ is a morphism in $\mathcal M$
and $g\co F(f)(v)\to w$ is a morphism in $F(n)$.) \qed

\medskip
By combining lemma~\ref{lem-longawait} with the locality result of \cite[\S4.2]{WeissSen1} for $\holim~\Theta_1$~,
we obtain a similar locality result for $\holim~\Theta$.
\end{appendices}

\bibliographystyle{amsplain}

\begin{thebibliography}{10}
\bibitem{Andrade} Ricardo Andrade, Ph.D thesis, MIT 2010.
\bibitem{Adams} J.~F.~Adams, \emph{Infinite loop spaces}, Annals of Math. Studies 90, Princeton Univ. Press 1978.
\bibitem{BoavidaAlpine} P.~Boavida de Brito, \emph{Segal objects and the Grothendieck construction}, in: An Alpine
bouquet of algebraic topology, 19--44, Contemp. Math. 708, Amer.Math.Soc., Providence RI, 2018.
\bibitem{BoavidaWeissLong} P.~Boavida de Brito and M.S.~Weiss, \emph{Spaces of smooth embeddings and configuration
categories}, J.Topol. 11 (2018), 65--143.
\bibitem{DwyKa2}  W.~G.~Dwyer and D.~Kan, \emph{Function complexes in homotopical algebra}, Topology 19 (1980), 427--440.
\bibitem{Hirschhorn} P.~S.~Hirschhorn, \emph{Model categories and their localizations}, Math. Surveys and Monographs vol.~99,
Amer.Math.Soc., 2002.
\bibitem{Hovey} M.~Hovey, \emph{Model categories},
Mathematical Surveys and Monographs, 63. Amer.Math.Soc., Providence, RI, 1999. xii+209 pp.
\bibitem{Rezk} C.~Rezk, \emph{A model for the homotopy theory of homotopy theory}, Trans.Amer.Math.Soc. 353 (2001), 973-1007.
\bibitem{Segal} G.~Segal, \emph{Categories and cohomology theories}, Topology 13 (1974), 293--312.
\bibitem{Thomason} R.W.~Thomason, \emph{Homotopy colimits in the category of small categories}, Math. Proc. Cambridge Philos. Soc. 85 (1979),
91–-109
\bibitem{TillSimp} S.~Tillmann, \emph{Occupants in simplicial complexes}, arXiv:1711.07107
\bibitem{WeissSen1} S.~Tillmann and M.~S.~Weiss, \emph{Occupants in manifolds}, in: Manifolds and $K$-theory,
237--259, Contemp.Math. 682, Amer.Math.Soc., Providence RI, 2017.
\bibitem{WeissDalian} M.~Weiss, \emph{Dalian notes on Pontryagin classes}, arXiv:1507.00153
\end{thebibliography}

\end{document}